\documentclass[11pt]{article}

\usepackage[utf8]{inputenc}
\usepackage{geometry}

\usepackage{amsmath,amssymb,amsfonts,amscd,amsthm,mathrsfs,mathtools,bbm,stmaryrd,mathabx,euscript,latexsym}

\usepackage{xcolor}
\usepackage{tikz-cd}
\usepackage{braket}
\usepackage{aliascnt}
\usepackage{thmtools}
\usepackage{ulem}
\usepackage{colortbl}
\usepackage{url}
\usepackage{authblk}
\usepackage{enumerate}
\usepackage[backref=page]{hyperref}
\usepackage{tikz}
\usepackage{graphicx}
\usepackage{todonotes}

\definecolor{lavender}{rgb}{0.9, 0.9, 0.98}
\definecolor{UGARed}{rgb}{0.729,0.047,0.184}

\hypersetup{
	colorlinks,
	linkcolor=UGARed,
	citecolor=blue,
	pdftitle={},
	bookmarksopen=true,
	bookmarksnumbered=true
}

\usepackage{cleveref}
\creflabelformat{enumi}{#2(#1)#3}
\crefformat{enumi}{(#2#1#3)}
\crefformat{equation}{(#2#1#3)}
\Crefformat{equation}{(#2#1#3)}
\crefformat{inequality}{(#2#1#3)}
\Crefformat{Inequality}{(#2#1#3)}

\newcommand\nnfootnote[1]{%
	\begin{NoHyper}
		\renewcommand\thefootnote{}\footnote{#1}%
		\addtocounter{footnote}{-1}%
	\end{NoHyper}
}

\definecolor{UGARed}{rgb}{0.729,0.047,0.184}

\theoremstyle{plain}
\newtheorem{thm}{Theorem}[section]

\newtheorem{lem}[thm]{Lemma}

\newtheorem{prop}[thm]{Proposition}

\theoremstyle{definition}
\newtheorem{defn}[thm]{Definition}

\newtheorem{rmk}[thm]{Remark}

\theoremstyle{plain}

\newtheorem{corollary}[thm]{Corollary}
\newtheorem{lemma}[thm]{Lemma}

\newtheorem{proposition}[thm]{Proposition}
\newtheorem{remark}[thm]{Remark}

\newtheorem{counter example}[thm]{Counter Example}
\newtheorem{question}[thm]{Question}

\newtheorem*{notation}{Notation}

\numberwithin{equation}{section}

\theoremstyle{plain}

\DeclareMathOperator\CC{{\mathbb C}}

\DeclareMathOperator\NN{{\mathbb N}}

\DeclareMathOperator\ZZ{{\mathbb Z}}

\DeclareMathOperator\Bl{{\mathcal B}}

\DeclareMathOperator\Hl{{\mathcal H}}

\DeclareMathOperator\Kl{{\mathcal K}}

\DeclareMathOperator\Rl{{\mathcal R}}

\DeclareMathOperator\Ul{{\mathcal U}}

\DeclareMathOperator{\Mor}{Mor}
\DeclareMathOperator{\id}{id}

\DeclareMathOperator{\Perm}{Perm}
\DeclareMathOperator{\PERM}{PERM}
\DeclareMathOperator{\Homeo}{Homeo}
\DeclareMathOperator{\HOMEO}{HOMEO}
\DeclareMathOperator{\Deck}{Deck}
\DeclareMathOperator{\DECK}{DECK}

\DeclareMathOperator{\Sym}{Sym}
\DeclareMathOperator{\Qut}{Qut}

\title{\bf Quantum Fundamental Groups of Spaces 
\\
with Finite-Sheeted Universal Covers}
\author{Shuzhou Wang\footnote{Department of Mathematics, University of Georgia, Athens, GA 30602, USA \href{mailto:szwang@uga.edu}{szwang@uga.edu}}, 
	Zhenhua Wang\footnote{Department of Mathematics, University of Georgia, Athens, GA 30602, USA \href{mailto:ezrawang@uga.edu}{ezrawang@uga.edu}}, and 
	Matthew Zediker\footnote{Department of Mathematics, University of Georgia, Athens, GA 30602, USA \href{mailto:Matthew.Zediker@uga.edu}{Matthew.Zediker@uga.edu}}}

\begin{document}
\maketitle

\begin{abstract}
 Discrete quantum groups have attracted much interest, but there was a lack of natural examples until Voigt constructed quantum permutations of finite and infinite sets in \cite{Voigt23}. In this paper, we construct the quantum fundamental groups that contain the classical fundamental groups as subgroups. This is a new class of natural examples of discrete quantum groups arising from symmetries of covering spaces in topology. We characterize quantum fundamental groups for a wide class of examples, including $\mathbb{RP}^2$. 
 We show that these quantum fundamental groups are genuinely quantum for any space 
 with finite fundamental groups that are of non prime order. For any odd prime order $p$, we construct a space with fundamental group of order $p$ whose quantum fundamental group is genuinely quantum.
\end{abstract}

\nnfootnote{\textit{Mathematics Subject Classification \rm{(2020)}:} {Primary 46L67; Secondary 46M15, 55R10, 57M05.}}

\nnfootnote{\textit{Keywords:}  Discrete quantum groups, quantum fundamental group, quantum deck transformation.}

\section{Introduction}

The theory of quantum groups emerged in the 1980s from two 
directions. On the one hand, Drinfeld and Jimbo \cite{Dr86,Jimb85} introduced quantum deformations of universal enveloping algebras of Lie algebras. On the other hand, Woronowicz \cite{Wor87, Wor87-SU2, Wor88}  developed the axioms of compact matrix quantum groups, and constructed quantum deformations of Lie groups $SU(n)$ as important examples.

Subsequently, through the work of the first author, compact matrix quantum groups that are not deformations of classical Lie groups were found from the study of intrinsic quantum symmetries. In particular, Wang introduced quantum permutation groups, which can be viewed as quantum automorphism groups of finite sets, or equivalently of finite-dimensional commutative \(C^*\)-algebras \cite{Wan98}. This opened a broader perspective in which quantum groups arise naturally from quantum symmetries of mathematical structures, rather than solely as deformations of classical Lie groups.

The theory of discrete quantum groups was developed in duality with the theory of compact quantum groups. Discrete quantum groups were originally defined as duals of compact quantum groups in Podle\'s and Woronowicz 
 \cite{PoWo90}, while intrinsic axiomatic formulations of discrete quantum groups were developed by Effros and Ruan \cite{EffrosRuan94} and Van Daele \cite{Vdl-discr96}. 
For a long period, concrete examples of discrete quantum groups were studied as duals of compact quantum groups without intrinsic independent constructions. 

An important new direction was initiated by Voigt \cite{Voigt23} 
through his introduction of infinite quantum permutation groups, extending Wang's quantum permutation groups from finite sets to infinite discrete sets. This provides a direct quantum-symmetry approach to construction of examples of discrete quantum groups, rather than starting with a compact quantum group and subsequently passing to its dual.   
In \cite{Voigt26},  Voigt further introduced quantum deck transformations using the same method. 

Our work follows this same direction, building on  \cite{Voigt23, Voigt26}. 
We start from a geometric symmetry problem, namely quantum symmetries of a covering space  
$
p:\widetilde{X}\longrightarrow X
$
of locally compact Hausdorff spaces that preserve its fibers. More precisely, we consider 
finite-dimensional quantum homeomorphisms $(\Hl_\sigma, (u^\sigma_{\tilde{x},\tilde{y}}))$
of the total space 
$\widetilde{X}$ 
satisfying the fiber-preserving condition

$$
u^\sigma_{\tilde{x},\tilde{y}}=0
\qquad\text{whenever}
\qquad
p(\tilde{x})\neq p(\tilde{y}).
$$

These objects are called 
{\bf finite-dimensional quantum deck transformations} of the covering 
and form a full subcategory of the category of finite-dimensional quantum permutations of total space $\widetilde{X}$ first introduced in \cite{Voigt23}.    
This is a \(C^*\)-tensor category with conjugates, 
called a rigid \(C^*\)-tensor category in \cite{Nesh_Tuset13,Voigt23,Voigt26}.  
Therefore Tannaka--Krein reconstruction can be invoked to produce a 
discrete quantum group  $\Qut^+_\delta(\widetilde{X}|X)$, called the 
\textbf{quantum deck transformation group} of the covering space 
$
p:\widetilde{X}\longrightarrow X
$. 
When the covering space is universal for a space $X$ satisfying appropriate conditions, 
we call the associated {quantum deck transformation group} 
the \textbf{quantum fundamental group} of $X$, and denote it by $\pi_1^q(X)$.  
We first describe $\pi_1^q(\mathbb{RP}^2)$. 
Then we extend the result to describe $\pi_1^q(X)$ for an appropriate class of spaces $X$ with finite-sheeted 
universal covers, and investigate its general properties.  

An outline of the paper is as follows. 
In \Cref{sec:prel}, we recall basic notions such as quantum permutations of a set \cite{Wan98, Voigt23}, quantum homeomorphisms, and quantum deck transformation group of a 
covering \cite{Voigt26}. 
Using the latter, in \Cref{sec:qfg}, 
we define the {quantum fundamental group} 
$\pi_1^q(X)$ of a space $X$ satisfying appropriate conditons.  
In \Cref{sec:qfgrp2}, we determine the quantum fundamental group $\pi_1^q(\mathbb{RP}^2)$ 
by describing the set of equivalence classes of 
irreducible finite-dimensional quantum deck transformations. 
This is an uncountable set and therefore 
$\pi_1^q(\mathbb{RP}^2)$  is an infinite discrete quantum group. 
As 
$\pi_1(\mathbb{RP}^2)$ is the two element group $\ZZ / 2 \ZZ$,  
the infinite discrete quantum group 
$\pi_1^q(\mathbb{RP}^2)$  can be viewed as a quantum version of 
$\ZZ / 2 \ZZ$.  
\Cref{sec:qfgrp2} serves as an illuminating guide to the determination of 
the more general quantum fundamental group $\pi_1^q(X)$ in \Cref{sec:qfgfscsp}  
using a similar method,  where $X$ is a path-connected, 
locally path-connected and semi-locally simply connected 
compact Hausdorff space with a finite fundamental group $\pi_1(X)$. 
Results in these two sections show that the discrete quantum group $\pi_1^q(X)$ is not a group and 
it properly contains the ordinary fundamental group $\pi_1(X)$.  
In the light of results of Sections \ref{sec:qfgrp2} and \ref{sec:qfgfscsp}, 
we investigate conditions when $\pi_1^q(X)$ is genuinely quantum, 
 i.e., when $\pi_1^q(X)$ is not the dual of a 
compact group, in \Cref{sec:clas}. 
 In particular, we show that $\pi_1^q(X)$ is genuinely quantum if   
 $\pi_1(X)$ is of non prime order; and 
 for any odd prime order $p$, 
 we construct a space with fundamental group of order $p$ whose quantum fundamental group is genuinely quantum. 
In \Cref{sec:galois}, we discuss a possible quantum Galois correspondence theory and open problems.

\section{Preliminaries}\label{sec:prel}

\noindent 
{\bf I. Quantum Permutations and Quantum Homeomorphisms}

\medskip

First, we recall the compact quantum permutation groups introduced in \cite{Wan98}. 
For basic notions on compact quantum groups, 
we refer the reader to \cite{Wor87,Baaj-Skan93,Wor98}. 

\begin{defn}[Quantum Permutation Group of a Finite Set \cite{Wan98}]\label{cqperm}
	The compact quantum group $S_N^+$ is the universal compact quantum permutation group on $N$ points. 
    The algebra of functions on $S_N^+$ is the universal (unital) $C^*$-algebra 
	$C(S_N^+)$ with generators $u_{ij}$, $1 \le i,j  \le  N$,  and relations  
	\begin{eqnarray}
& &	u_{ij}=u_{ij}^*=u_{ij}^2, \label{MagicProj}
\\
& &
\sum_{k=1}^Nu_{kj}=1=\sum_{k=1}^Nu_{ik}, \label{MagicProjSum}
	\end{eqnarray}
	and the comultiplication $\Delta: C(S_N^+) \to C(S_N^+) \otimes C(S_N^+)$ is given  
	by \[\Delta(u_{ij})=\sum_{k=1}^Nu_{ik}\otimes u_{kj}.\] 
Here and throughout this paper, we use spatial $C^*$-tensor product. 	
\end{defn}

Voigt extended the above construction to define the discrete quantum permutation group of any set \cite{Voigt23}.

\begin{defn}[Quantum Permutation of a Set \cite{Voigt23}]
Let $\widetilde{X}$ be a set. A \textbf{quantum permutation} of $\widetilde{X}$ is a
pair $\sigma=(\Hl_\sigma,u^\sigma)$ consisting of a Hilbert space
$\Hl_\sigma$ and a family of self-adjoint projections
\[
u^\sigma=(u^\sigma_{\tilde{x},\tilde{y}})_{\tilde{x},\tilde{y}\in\widetilde{X}},
\qquad
u^\sigma_{\tilde{x},\tilde{y}}\in\Bl(\Hl_\sigma),
\]
such that for all $\tilde{x},\tilde{y}\in\widetilde{X}$,
\[
\sum_{\tilde{y}\in\widetilde{X}}u^\sigma_{\tilde{x},\tilde{y}}
=
I_{\Hl_\sigma}
=
\sum_{\tilde{x}\in\widetilde{X}}u^\sigma_{\tilde{x},\tilde{y}},
\]
where the sums converge in the strong operator topology (SOT) of
$\Bl(\Hl_\sigma)$. 
Here we adapt the notation in \cite{Voigt26} for $u^\sigma_{\tilde{x},\tilde{y}}$.  
Note that in the notation of \cite{Voigt23}, 
$u^\sigma_{\tilde{x},\tilde{y}} = u^\sigma_{\tilde{x} \tilde{y}}$.

The \textbf{quantum permutation $\sigma$ is called finite-dimensional} if
$\Hl_\sigma$ is finite-dimensional.
\end{defn}

One may define direct sums, subobjects, tensor products, and conjugates of quantum permutations; 
see \cite{Voigt23} for details.

\begin{defn}[Category of 
Finite-Dimensional 
Quantum Permutations \cite{Voigt23}]
The objects of the category 
$
\displaystyle {\PERM}^+(\widetilde{X})
$
of quantum permutations
are 
$
\sigma=(\Hl_\sigma,u^\sigma)
$ of $\widetilde{X} $ as defined above.

For two objects
$
\sigma=(\Hl_\sigma,u^\sigma),
\; 
\tau=(\Hl_\tau,u^\tau),$
a morphism in $\Mor(\sigma, \tau)$ 
is a bounded operator
$
T\in\mathcal{B}(\Hl_\sigma,\Hl_\tau)
$
such that
\[
T u^\sigma_{\tilde{x},\tilde{y}}
=
u^\tau_{\tilde{x},\tilde{y}}T,
\qquad
\tilde{x},\tilde{y}\in\widetilde{X}.
\]

The \textbf{category of finite-dimensional quantum permutations} 
$
{\Perm}^+(\widetilde{X})
$
is the full $C^*$-subcategory of 	
$
{\PERM}^+(\widetilde{X})
$
with finite dimensional quantum permutations as the objects.
\end{defn}

Because of restriction to finite-dimensional objects,  
$
{\Perm}^+(\widetilde{X})
$
is a rigid $C^*$-tensor category, and Tannaka-Krein reconstruction theorem 
\cite{Wor88, Nesh_Tuset13, Vdl-discr96, Vdl-duality98, Vdl-discr25} 
can be applied to obtain 	

\begin{defn}[Quantum Permutation Group of an Infinite Set \cite{Voigt23}]\label{qpermgp}
The  \textbf{quantum permutation group} $\Sym^+(\widetilde{X})$ of $\widetilde{X}$ is the discrete quantum group 
associated to 
$
{\Perm}^+(\widetilde{X})
$
via 
Tannaka-Krein reconstruction. 	
\end{defn}

Let us recall \cite{Voigt23} the quantum algebra of functions $C_0(\Sym^+(\widetilde{X}))$  on $\Sym^+(\widetilde{X})$ and the coproduct 
$\Delta: C_0(\Sym^+(\widetilde{X})) \to 
M(C_0(\Sym^+(\widetilde{X})) \otimes C_0(\Sym^+(\widetilde{X})))$ as a nondegenerate 
$C^*$-algebra morphism. 
Let 
${\rm Irred}(\widetilde{X})$ or (if no confusion arises) simply ${\rm Irred}$ be a complete set of representatives $\sigma$ of unitary equivalence classes of irreducible objects in 
$\Perm^+(\widetilde{X})$. 
Then  
\begin{align*}
C_0(\Sym^+(\widetilde{X}))& := \bigoplus_{\sigma \in {\rm Irred} } B(H_{\sigma}),
\\
 \Delta (u_{\tilde{x}, \tilde{y}})& := \sum_{\tilde{z} \in \widetilde{X}} 
u_{\tilde{x}, \tilde{z}} \otimes u_{\tilde{z}, \tilde{y}},     
\end{align*}
where $u_{\tilde{x}, \tilde{y}} \in M(C_0(\Sym^+(\widetilde{X})))$ with components 
$u_{\tilde{x}, \tilde{y}}^\sigma$, 
i.e., $u = (u_{\tilde{x}, \tilde{y}})_{\tilde{x}, \tilde{y} \in \widetilde{X}}$ is the universal quantum permutation.

\begin{rmk}
The novel approach of Voigt \cite{Voigt23} recalled here that is different from that in 
Wang \cite{Wan98} is that Voigt considered the representations of the 
{\bf quantum function algebra} of the discrete quantum group to be defined using the 
category of quantum permuations, 
whereas Wang considered representations of 
the {\bf compact quantum automorphism group} to be defined via quantum transformation groups. 
\end{rmk}

\begin{rmk}
The discrete quantum groups $\Sym^+(\widetilde{X})$ are examples of locally compact quantum groups in the sense of Kustermans and Vaes \cite{KusVaes00, KusVaes03}, though  
the existence of left and right Haar measures on discrete quantum groups, which are postulated in \cite{KusVaes00},  
are a consequence of the axioms of discrete quantum groups 
\cite{Vdl-discr96, Vdl-discr25}
just as in the case of compact quantum groups \cite{Wor87,Baaj-Skan93, Vdl-Haar,Wor98}. 
\end{rmk}

We recall the fundamental notions of quantum subgroups of a discrete quantum group and discrete quantum group morphisms, which has no ambiguity because of the results in 
Daws-Kasprzak-Skalski-So{\l}tan \cite{Daws12} in our setting based on 
Vaes \cite{Vaes05} and Meyer-Roy-Woronowicz \cite{Meyer12}. 
See also \cite{Wan95, Wan13} for the case of compact quantum groups.

\begin{defn}
    Let $(C_0(\Gamma_1),\Delta_1)$ and $(C_0(\Gamma_2),\Delta_2)$ denote discrete quantum groups. A \textbf{morphism of discrete quantum groups} $\iota:\Gamma_1\to\Gamma_2$ is a $*$-homomorphism $\iota^*:C_0(\Gamma_2)\to M(C_0(\Gamma_1))$ which intertwines the comultiplications: \[(\iota^*\otimes\iota^*)\circ\Delta_2=\Delta_1\circ\iota^*.\] We say $\Gamma_1$ is a \textbf{subgroup} of $\Gamma_2$ and write $\Gamma_1\leq\Gamma_2$ if there exists a surjective $*$-homomorphism $\iota^*:C_0(\Gamma_2)\twoheadrightarrow C_0(\Gamma_1)$ which intertwines the comultiplications as above.
\end{defn}

\begin{rmk}
Note that the one-dimensional quantum permutations correspond to classical set of bijections $\Sym(\widetilde{X})$, so we always have the containment $\Sym(\widetilde{X})\leq\Sym^+(\widetilde{X})$.
\end{rmk}

\begin{defn}[Quantum Homeomorphisms \cite{Voigt26}]\label{qhomeo}
Let $\widetilde{X}$ be a locally compact Hausdorff space and let
$\widetilde{X}_{\delta}$ denote the same set equipped with the discrete
topology. A quantum permutation 
$
\sigma=(\Hl_\sigma,u_\sigma)
$
of 
$\widetilde{X}_{\delta}$ is called a 
\textbf{quantum homeomorphism of $X$} if 
for every $f\in C_0(\widetilde{X})$ both maps
\[
\widetilde{X}\ni\tilde{y}
\longmapsto
\sum_{\tilde{x}\in\widetilde{X}}
f(\tilde{x})
u^\sigma_{\tilde{x},\tilde{y}}
\]
and
\[
\widetilde{X}\ni\tilde{x}
\longmapsto
\sum_{\tilde{y}\in\widetilde{X}}
f(\tilde{y})
u^\sigma_{\tilde{x},\tilde{y}}
\]
belong to
$
C_0(\widetilde{X},\Bl(\Hl_\sigma)).
$
\end{defn}

\begin{remark}
{\rm 
	The above is equivalent to the requirement that both $\sigma$ and $\bar{\sigma}$ are 
	continuous as defined in \cite{Voigt26}}.
	\end{remark}

	\begin{defn}[Category of Finite-Dimensional Quantum Homeomorphism \cite{Voigt26}]
The category 
\\
$
{\HOMEO}^+(\widetilde{X} )
$ 
is the full $C^*$-tensor subcatory of 
$
{\PERM}^+(\widetilde{X} )
$
with quantum homeomorphisms of $\widetilde{X}$ as objects. 

The category 
$
{\Homeo}^+(\widetilde{X} )
$ 
is the full 
$C^*$-tensor subcatory of 
$
{\HOMEO}^+(\widetilde{X} )
$
with finite dimensional quantum homeomorphisms of $\widetilde{X}$ as the objects. 
\end{defn}

\begin{rmk}
As $\Perm^+(\widetilde{X}_\delta)$\footnote{By a slight abuse of notation, one may simply use 
	$\Perm^+(\widetilde{X})$ for $\Perm^+(\widetilde{X_\delta})$ since $\Perm^+$ 
	is only defined for discrete spaces.}, 
$
{\Homeo}^+(\widetilde{X})
$ 
is a rigid $C^*$-tensor category, and a full subcatory of 
$
{\PERM}^+(\widetilde{X})
$.
\end{rmk}

\begin{defn}[Quantum Homeomorphism Group \cite{Voigt26}]
The  \textbf{quantum homeomorphism group} $\Qut^+_\delta(\widetilde{X})$ 
of $\widetilde{X}$ is the discrete quantum group 
associated to the rigid $C^*$-tensor category 
$
{\Homeo}^+(\widetilde{X})
$
via 
Tannaka-Krein reconstruction. 	
\end{defn}

\noindent
{\bf II. Quantum Deck Transformation Group}

\begin{defn}[Quantum Deck Transformation \cite{Voigt26}]
Let $X$ be a locally compact Hausdorff space and 
$
p:\widetilde{X}\longrightarrow X
$
be a covering map. 
For ease of notation, denote the covering by 
$\widetilde{X} | X.$

A 
\textbf{quantum deck transformation}
of the covering
$\widetilde{X} | X$  
is a quantum homeomorphism 
\\
$
\sigma=(\Hl_\sigma,u^\sigma)
$
of the total space $\widetilde{X}$ such that
\[
u^\sigma_{\tilde{x},\tilde{y}}=0
\quad
\text{whenever }
p(\tilde{x})\neq p(\tilde{y}),
\; \; \text{where}
\; \; 	
\tilde{x},\tilde{y}\in\widetilde{X}.
\]

Equivalently, the associated quantum permutation preserves each fiber 
\[
p^{-1}(x),\qquad x\in X,
\]
and hence only permutes points lying in the same fiber of the covering.

$\sigma$ is called a {\bf finite-dimensional quantum deck transformation} if
$\Hl_\sigma$ is finite-dimensional.

\end{defn}

\begin{defn}[Category of Finite-Dimensional Quantum Deck Transformations \cite{Voigt26}]
The objects of the category
$
{\DECK}^+(\widetilde{X} | X)
$
are quantum deck transformations
$
\sigma=(\Hl_\sigma,u_\sigma)
$
of the covering 
$\widetilde{X} | X$.

For two 
deck transformations
$
\sigma=(\Hl_\sigma,u^\sigma),
\; 
\tau=(\Hl_\tau,u^\tau),$
the space of morphism $\Mor(\sigma, \tau)$ 
is the same as that for $\PERM^+(\widetilde{X}_\delta)$ when 
$\sigma$ and 
$\tau$
are viewed as objects of $\PERM^+(\widetilde{X}_\delta)$\footnote{The only difference here from Definition 5.12 in \cite{Voigt26} is that $T$ is 
	independent of $\tilde{x}\in\widetilde{X}$ so that $\Mor(\sigma, \tau)$ is finite dimensional when $\sigma, \tau$ are finite dimensional and Tannaka-Krein reconstruction can be invoked.}.

The \textbf{category of finite-dimensional quantum deck transformations} 
$
{\Deck}^+(\widetilde{X} | X)
$
is the full $C^*$-tensor subcategory of 	
$
{\DECK}^+(\widetilde{X} | X)
$
with finite dimensional quantum deck transformations as the objects.
\end{defn}

\begin{rmk}
By definition,  
${\Deck}^+(\widetilde{X} | X)$
is a full rigid $C^*$-tensor subcategory of 
${\Homeo}^+(\widetilde{X})$ and ${\Perm}^+(\widetilde{X})$.	In particular, one-dimensional quantum deck transformations correspond to classical deck transformations.
\end{rmk}

\begin{defn}[Quantum Deck Transformation Group \cite{Voigt26}]
The  \textbf{quantum deck transformation group} $\Qut^+_\delta(\widetilde{X}|X)$ 
of $\widetilde{X}|X$ is the discrete quantum group 
associated to the rigid $C^*$-tensor category 
${\Deck}^+(\widetilde{X} | X)$
via 
Tannaka-Krein reconstruction. 	
\end{defn}

\section{Quantum Fundamental Groups}\label{sec:qfg}
\noindent

Henceforth, $X$ is a locally compact Hausdorff space which is path-connected, locally path-connected, and semi-locally simply connected. Classically, the fundamental group of $X$ is isomorphic to the deck transformation group of its universal covering space. 
It is natural to define a quantum fundamental group from these quantum deck transformations. 
{\it A priori}, it is not obvious that such a definition would be independent of the 
choice of universal cover. We show this is indeed the case below. 

Let $p_1:\widetilde{X}_1\rightarrow X$ and $p_2:\widetilde{X}_2\rightarrow X$ be covering maps. Recall the covering spaces $\widetilde{X}_1|X$ and $\widetilde{X}_2|X$ are isomorphic if there is a homeomorphism $F:\widetilde{X}_1\rightarrow\widetilde{X}_2$ such that $p_2\circ F=p_1$.

\begin{proposition}
\label{welldefined}
Let $X$ be a path-connected, locally path-connected, and semi-locally simply-connected locally compact Hausdorff space. Let $p_1:\widetilde{X}_1\rightarrow X$ and $p_2:\widetilde{X}_2\rightarrow X$ be covering maps. Suppose there is a covering space isomorphism $F:\widetilde{X}_1\rightarrow\widetilde{X}_2$. Then, $\Qut^+_\delta(\widetilde{X}_1|X)\cong\Qut^+_\delta(\widetilde{X}_2|X)$ as discrete quantum groups. 
\end{proposition}

\begin{proof}
Let $\sigma=(\Hl_\sigma,u^\sigma)\in\text{Deck}^+(\widetilde{X}_1)$ be an object. We define the pair $F_*\sigma:=(\Hl_\sigma,v^{F_*\sigma})$, where 
\[
v^{F_*\sigma}_{\tilde{x}_2,\tilde{y}_2}:=u^\sigma_{F^{-1}(\tilde{x}_2),F^{-1}(\tilde{y}_2)}\quad \mbox {for each}\,\, \tilde{x}_2,\, \tilde{y}_2\in\widetilde{X}_2. 
\]

Since $F^{-1}$ is continuous, and $\sigma$ is a quantum homeomorphism, $F_*\sigma$ is a quantum homeomorphism. If $p_2(\tilde{x}_2)\neq p_2(\tilde{y}_2)$, then $p_1(F^{-1}(\tilde{x}_2))\neq p_1(F^{-1}(\tilde{y}_2))$ by the condition $p_1=p_2\circ F$. Since $\sigma$ is a quantum deck transformation, $v_{\tilde{x}_2,\tilde{y}_2}^{F_*\sigma}=u_{F^{-1}(\tilde{x}_2),F^{-1}(\tilde{y}_2)}^\sigma=0$. Thus, $F_*:\sigma\mapsto F_*\sigma$ defines a map from the objects of $\text{Deck}^+(\widetilde{X}_1)$ to the objects of $\text{Deck}^+(\widetilde{X}_2)$. By defining $F_*$ to be the identity on morphisms, 
\[
\text{Mor}(F_*(\sigma),F_*(\tau))=\text{Mor}(\sigma,\tau),
\]
we obtain a functor $F_*$.

Then, $F_*$ induces a one-to-one correspondence from equivalence classes of irreducible objects in $\text{Deck}^+(\widetilde{X}_1)$ to equivalence classes of irreducible objects in $\text{Deck}^+(\widetilde{X}_2)$. Bijectivity of $F$ ensures compatibility of $F_*$ with the tensor products in the respective categories. 

Now, we can decompose the universal matrix entries $v_{\tilde{x},\tilde{y}}$ of $\Qut^+_\delta(\widetilde{X}_2|X)$ as 
\[
v_{\tilde{x},\tilde{y}}=\sum_{\sigma}v^{
F_*\sigma}_{\tilde{x},\tilde{y}},
\]
where $\sigma$ indexes a complete set of representatives of unitary equivalence classes of irreducible objects in the category $\text{Deck}^+(\widetilde{X}_1|X)$ (with convergence in the SOT). 

Recall $C_0(\Qut^+_\delta(\widetilde{X}_1|X))$ and $C_0(\Qut^+_\delta(\widetilde{X}_2|X))$ are direct sums of matrix algebras with respective matrix units $e_{ij}^\sigma$ and ${f}_{kl}^\tau$. Define a $*$-homomorphism $\phi^*:C_0(\Qut^+_\delta(\widetilde{X}_1|X))\to C_0(\Qut^+_\delta(\widetilde{X}_2|X))$ by 
\[
\phi^*(e_{ij}^\sigma)={f}_{ij}^{F_*\sigma}.
\] 
This satisfies $\phi^*(u_{\tilde{x},\tilde{y}})=v_{\tilde{x},\tilde{y}}$, where 
\[
u_{\tilde{x},\tilde{y}}=\sum_{\sigma}u^{\sigma}_{\tilde{x},\tilde{y}}.\]
Since $F_*$ preserves tensors and irreducibility, $\phi^*$ is compatible with the comultiplications, given by 
\[\Delta(u_{\tilde{x},\tilde{y}})=\sum_{\tilde{z} \in \widetilde{X}_1} 
u_{\tilde{x}, \tilde{z}} \otimes u_{\tilde{z}, \tilde{y}}
\text{\quad and \quad}\Delta(v_{\tilde{x},\tilde{y}})=\sum_{\tilde{z} \in \widetilde{X}_2} 
v_{\tilde{x}, \tilde{z}} \otimes v_{\tilde{z}, \tilde{y}}.
\] 
Thus, $\phi^*$ defines an isomorphism $\phi$ of discrete quantum groups $\Qut^+_\delta(\widetilde{X}_2|X)\cong\Qut^+_\delta(\widetilde{X}_1|X)$.
\end{proof}

This proposition justifies the following definition.

\begin{defn}[Quantum Fundamental Group]
We call $\Qut^+_\delta(\widetilde{X}|X)$ the \textbf{quantum fundamental group} of $X$ and denote it by $\pi_1^q(X)$.
\end{defn}

In the context of $\pi_1^q(X)$, similar to \Cref{qpermgp}, the quantum algebra of functions $C_0(\pi_1^q(X))$  on $\pi_1^q(X)$ and the coproduct $\Delta$ are defined as follows.

Let 
${\rm Irred}_1^q(X)$ or (if no confusion arises) simply ${\rm Irred}_1^q$ be a complete set of representatives $\sigma$ of unitary equivalence classes of irreducible objects in ${\Deck}^+(\widetilde{X} | X)$. Then  

\[
C_0(\pi_1^q(X)) : = \bigoplus_{\sigma \in {\rm Irred}_1^q } B(H_{\sigma}),
\]
\[
\Delta (u_{\tilde{x}, \tilde{y}}) = \sum_{\tilde{z} \in \widetilde{X}} 
u_{\tilde{x}, \tilde{z}} \otimes u_{\tilde{z}, \tilde{y}}, 
\]
where $u_{\tilde{x}, \tilde{y}} \in M(C_0(\pi_1^q(X)))$ with components 
$u_{\tilde{x}, \tilde{y}}^\sigma$, 
i.e., $u = (u_{\tilde{x}, \tilde{y}})_{\tilde{x}, \tilde{y} \in \widetilde{X}}$ is the universal quantum deck transformation.

We show this discrete quantum group satisfies a universal property. First, we establish the following notions. See also \cite{Voigt23, Voigt26}.

\begin{defn}\label{actions}
Let $\Gamma$ be a discrete quantum group with a comultiplication $\Delta_\Gamma$ and $\widetilde{Y}$ be a universal cover for a space $Y$ with covering map $p:\widetilde{Y}\rightarrow Y$. 

\begin{enumerate}[\rm(i)]

\item An \textbf{action} of $\Gamma$ on $\widetilde{Y}$ is an injective morphism of $C^*$-algebras $\alpha:C_0(\widetilde{Y})\rightarrow M(C_0(\Gamma)\otimes C_0(\widetilde{Y}))$ satisfying the following: 
\begin{enumerate}[\rm(i)]
    \item $(\Delta_\Gamma\otimes\id)\circ\alpha=(\id\otimes\:\alpha)\circ\alpha$, and
    \item $\alpha(C_0(\widetilde{Y}))(C_0(\Gamma)\otimes 1)$ is dense in $C_0(\Gamma)\otimes C_0(\widetilde{Y})$.
\end{enumerate}

\item An action $\alpha_\delta:C_0(\widetilde{Y}_\delta) \to M(C_0(\Gamma)\otimes C_0(\widetilde{Y}_\delta))$ is \textbf{by permutations} if the counting measure $\gamma$ on $C_0(\widetilde{Y}_\delta)$ is invariant in the sense that 
\[\gamma((\omega\otimes\id)(\alpha_\delta(f)))=\gamma(f)\omega(1)\] 
for all $f\in C_0(\widetilde{Y}_\delta)$ and $\omega$ a finitely-supported linear functional on $C_0(\Gamma)$. The universality in \cite{Voigt23} shows we can write $\alpha_\delta(e_{\tilde{x}})=:\displaystyle\sum_{\tilde{y}\in\widetilde{Y}}v_{\tilde{x},\tilde{y}}\,\otimes\,e_{\tilde{y}}$ for some uniquely determined $v_{\tilde{x},\tilde{y}}\in M(C_0(\Gamma))$ which satisfy the quantum permutation relations. 

\item \label{itm:acthmph}We say $\alpha:C_0(\widetilde{Y}) \to M(C_0(\Gamma)\otimes C_0(\widetilde{Y}))$ an action \textbf{by homeomorphisms} if there is an action by permutations $\alpha_\delta:C_0(\widetilde{Y}_\delta) \to M(C_0(\Gamma)\otimes C_0(\widetilde{Y}_\delta))$, called the \textbf{discretization} of $\alpha$, making the diagram \[\begin{tikzcd}
	C_0(\widetilde{Y}) \arrow{r}{\alpha} \arrow[swap]{d}{\iota^*} & M(C_0(\Gamma)\otimes C_0(\widetilde{Y}))\arrow{d}{\id\otimes\,\iota^*} \\ M(C_0(\widetilde{Y}_\delta)) \arrow[swap]{r}{\alpha_\delta} & M(C_0(\Gamma)\otimes C_0(\widetilde{Y}_\delta))
\end{tikzcd}\]
commute and 
\[\overline{\alpha}_\delta(e_{\tilde{y}}):=\sum_{\tilde{x}\in\widetilde{Y}}v_{\tilde{x},\tilde{y}}\otimes e_{\tilde{x}}\] 
defines an action $\overline{\alpha}:C_0(\widetilde{Y}) \to M(C_0(\Gamma)\otimes C_0(\widetilde{Y}))$ via 
\[(\id\otimes\:\iota^*)\circ\overline{\alpha}=\overline{\alpha}_\delta\circ\iota^*.\] 
Here and in the proof below, $\iota^*$ is induced by the (continuous) inclusion $\iota:\widetilde{X}_\delta\hookrightarrow\widetilde{X}$, i.e.,  
\[\iota^*(f)=\sum_{\tilde{y}\in\widetilde{Y}}f(\tilde{y})e_{\tilde{y}},\]
where $e_{\tilde{y}}\in C_0(\widetilde{Y}_\delta)$ is the characteristic function on $\tilde{y}\in\widetilde{Y}$, for all $f\in C_0(\widetilde{Y})$.

\item Finally, the action $\alpha$ of $\Gamma$ on $\widetilde{Y}$ is 
\textbf{by deck transformations} if it acts by homeomorphisms and the diagram 
\[
\begin{tikzcd}
	C_0(Y) \arrow{r}{\tau} \arrow[swap]{d}{p^*} & M(C_0(\Gamma)\otimes C_0(Y)) \arrow{d}{\id\otimes\:p^*} \\ M(C_0(\widetilde{Y})) \arrow[swap]{r}{\alpha}
	& M(C_0(\Gamma)\otimes C_0(\widetilde{Y}))
\end{tikzcd}
\] 
commutes, where $\tau$ is the trivial action defined by $\tau(f)=1\otimes f$ for $f\in C_0(Y)$. We call such an action $\alpha$ a \textbf{deck action} of $\Gamma$ on the covering space.
\end{enumerate}
\end{defn}

\begin{rmk}
    \Cref{actions} \cref{itm:acthmph} is used in \cite{Voigt26} to give an equivalent formulation of \Cref{qhomeo}.
\end{rmk}

\begin{prop}
\label{universality}
Let $X$ be a locally compact Hausdorff space which is path-connected, locally path-connected, and semilocally simply-connected. Fix a universal cover $p:\widetilde{X}\rightarrow X$. Then, $\pi_1^q(X)$ is the universal discrete quantum group acting on $\widetilde{X}$ by quantum deck transformations. 

That is, there is a deck action $\beta$ of $\pi_1^q(X)$ on $\widetilde{X}$ such that, if $\Gamma$ is any discrete quantum group with a deck action $\alpha:C_0(\widetilde{X})\rightarrow M(C_0(\Gamma)\otimes C_0(\widetilde{X}))$, there exists a unique morphism of discrete quantum groups $\phi^*:C_0(\pi_1^q(X))\rightarrow M(C_0(\Gamma))$ satisfying $(\phi^*\otimes\id)\circ\beta=\alpha$.
\end{prop}

This proof mirrors the proof of universality in \cite{Voigt23}. The main difference is that additional relations in our setting must be checked.

\begin{proof}
First, note that 
$\beta_\delta:C_0(\widetilde{X}_\delta)\rightarrow M(C_0(\pi_1^q(X))\otimes C_0(\widetilde{X}_\delta))$ 
given by 
$\beta_\delta(e_{\tilde{x}})
=\sum_{\tilde{y}\in\widetilde{X}}u_{\tilde{x},\tilde{y}}\otimes e_{\tilde{y}}$ 
defines the discretization of the desired deck action 
$\beta:C_0(\widetilde{X})\rightarrow M(C_0(\pi_1^q(X))\otimes C_0(\widetilde{X}))$ 
of $\pi_1^q(X)$ on $\widetilde{X}$, 
where $u_{\tilde{x},\tilde{y}}$ are the matrix entries of the universal quantum deck transformation.

Let $\Gamma$ be a discrete quantum group acting on $\widetilde{X}$ by 
deck transformations, and denote the action by 
$\alpha:C_0(\widetilde{X})\rightarrow M(C_0(\Gamma)\otimes C_0(\widetilde{X}))$. 
Let $\alpha_\delta$ denote the discretization of $\alpha$. 
Then, we can write 
\[\alpha_\delta(e_{\tilde{x}})
=:\sum_{\tilde{y}\in\widetilde{X}}v_{\tilde{x},\tilde{y}}\,\otimes\,e_{\tilde{y}}\]
for some uniquely determined $v_{\tilde{x},\tilde{y}}\in M(C_0(\Gamma))$. 
Moreover, since $\alpha_\delta$ is an action by permutations, 
components of $v_{\tilde{x},\tilde{y}}$ are projections which satisfy the relations to be a quantum permutation. The components of the entries $v_{\tilde{x},\tilde{y}}$ also satisfy the condition in \Cref{qhomeo}, since $\alpha$ is an action by homeomorphisms compatible with $\alpha_\delta$. 

Next we check that $v_{\tilde{x},\tilde{y}}=0$ if $p(\tilde{x})\neq p(\tilde{y})$. We observe that
\[\alpha_\delta\circ\iota^*\circ p^*=(\id\otimes\,\iota^*)\circ(\id\otimes\,p^*)\circ\tau\]
by stacking the commutative diagrams above. For all $f\in C_0(X)$, 
\begin{align*}
	\alpha_{\delta}\circ\iota^*\circ p^*(f)&=\alpha_\delta\circ\iota^*\circ f\circ p=\alpha_\delta\Bigg(\sum_{\tilde{x}\in\widetilde{X}_\delta}f(p(\tilde{x}))e_{\tilde{x}}\Bigg)
   \\ &=\sum_{\tilde{x}\in\widetilde{X}_\delta}f(p(\tilde{x}))\Bigg(\sum_{\tilde{y}\in\widetilde{X}_\delta}v_{\tilde{x},\tilde{y}}\otimes e_{\tilde{y}}\Bigg)
    \\
&=\sum_{\tilde{y}\in\widetilde{X}_\delta}\Bigg(\sum_{\tilde{x}\in\widetilde{X}_\delta}f(p(\tilde{x}))v_{\tilde{x},\tilde{y}}\Bigg)\otimes e_{\tilde{y}},
\end{align*}
and 
\begin{align*}
(\id\otimes\,\iota^*)\circ(\id\otimes\, p^*)\circ\tau(f)
&=(\id\otimes\,\iota^*)\circ(\id\otimes\, p^*)(1\otimes f)\\
&=1\otimes\Bigg(\sum_{\tilde{x}\in\widetilde{X}_\delta}f(p(\tilde{x}))e_{\tilde{x}}\Bigg)\\
&=\sum_{\tilde{x}\in\widetilde{X}_\delta}f(p(\tilde{x}))1\otimes e_{\tilde{x}}.
\end{align*} 
Hence, for any $f\in C_0(X)$ and $\tilde{y}\in\widetilde{X}$, 
\[
\sum_{\tilde{x}\in\widetilde{X}_\delta}f(p(\tilde{x}))v_{\tilde{x},\tilde{y}}=f(p(\tilde{y}))1.
\]

Fix some $\tilde{x},\tilde{y}\in\widetilde{X}$, and suppose 
$p(\tilde{x})\neq p(\tilde{y})$. 
Choose a function $f\in C_0(X)$ with $f(p(\tilde{x}))=1$ and $f(p(\tilde{y}))=0$. 
Then the above identity becomes, 
\[\sum_{\tilde{z}}f(p(\tilde{z}))v_{\tilde{z},\tilde{y}}=f(p(\tilde{y}))1=0.
\] 
Multiplying with $v_{\tilde{x},\tilde{y}}$ on both sides, 
since $ f(p(\tilde{z})) \displaystyle v_{\tilde{z},\tilde{y}}v_{\tilde{x},\tilde{y}}
= f(p(\tilde{z})) \delta_{\tilde{x},\tilde{z}}v_{\tilde{x},\tilde{y}}$ and $\displaystyle f(p(\tilde{x}))=1$, 
we obtain 
\[
v_{\tilde{x},\tilde{y}}=0.
\]

From the definitions in \cite{Vdl-discr96, Vdl-discr25}, 
the $C^*$-algebra of functions on $\Gamma$ can be written as
\[
C_0(\Gamma)=\bigoplus_{\gamma\in I}\Bl(\Kl_{\gamma})
\]
for some index set $I$. We write 
\[
v_{\tilde{x},\tilde{y}}=\sum_{\gamma'\in I^{'}}v^{\gamma'}_{\tilde{x},\tilde{y}}
\] 
for some smallest $I'\subseteq I$. Then, for each $\gamma'\in I'$, the pair $(\Kl_{\gamma'},(v^{\gamma'}_{\tilde{x},\tilde{y}}))$ is a quantum deck transformation of $\tilde{X}|X$. Thus, each such $\gamma'\in I'$ is an object in $\text{Deck}^+(\widetilde{X}|X)$. 

Since each $\gamma'$ is equivalent to a direct sum of irreducible subobjects $\sigma\in\text{Irred}_1^q$, and $u_{\tilde{x},\tilde{y}}$ is a sum indexed by irreducibles, there is a well-defined $*$-homomorphism $\phi^*:C_0(\pi_1^q(X))\rightarrow M(C_0(\Gamma))$ by assigning the irreducible blocks in \[
C_0(\pi_1^q(X))=\bigoplus_{\sigma\in\text{Irred}_1^q}\Bl(\Hl_\sigma)
\] to the blocks of $C_0(\Gamma)$ corresponding to the decomposition of each $\gamma'\in I'$ into irreducibles, perhaps with multiplicities; 
for those $\sigma\in\text{Irred}_1^q$ not contained in any $\gamma'$, map the block $\Bl(\Hl_\sigma)$ to $0$. This satisfies the formula 
\[\phi^*(u_{\tilde{x},\tilde{y}})=v_{\tilde{x},\tilde{y}}.
\]
This construction is compatible with the comultiplications, since the $\gamma'$ are also objects in the rigid $C^*$-tensor category associated to $\Gamma$. Thus, we obtain a morphism of discrete quantum groups $\phi:\Gamma\rightarrow\pi_1^q(X)$. Uniqueness of $\phi$ follows from the uniqueness of the $v_{\tilde{x},\tilde{y}}$. 
\end{proof}

\begin{rmk}
    We can now say $\pi_1^q(X)$ is \textbf{the} discrete quantum deck transformation group of $\widetilde{X}|X$, as now any discrete quantum deck transformation group action on $\widetilde{X}$ must factor through the action $\beta$ defined above.
\end{rmk}

\begin{rmk}
Presently, the following are three prominent approaches to 
construct quantum groups from the functions on quantum group point of view, 
though they are not mutually exclusive:

(1) Define quantum function algebra on a quantum group by generators and relations similar to defining algebraic varieties using polynomial identities. Quantum subgroups are 
defined by adding more relations. Examples of this approach include 
$C(U_n^+)$, $C(O_n^+)$, $A_u(Q)$, $B_u(Q)$ in \cite{Wan95,VanWan96,Wan02b}.

(2) Define quantum function algebra 
on a quantum group by using representations of the quantum group to be defined. 
Quantum subgroups are 
defined by looking at specialized representations of the quantum group. 
Examples of this approach include \cite{Wan98,Bich03a,Ban05a,Ban05b,Ban-Sp09} 
and the vast literature on quantum automorphism groups of finite graphs.

(3) Define quantum function algebra on a quantum group by using the category of the  representations of the former. Quantum subgroups are 
defined using full subcategories. 
Examples of this approach include \cite{Voigt23,Voigt26} 
as recalled in Section 2 and the construction of 
quantum fundamental group $\pi_1^q(X)$ defined above. 
\end{rmk}

\section{The Quantum Fundamental Group of \texorpdfstring{$\mathbb{RP}^2$}{}} \label{sec:qfgrp2}

Here, we consider the basic example of $\mathbb{RP}^2$. We describe the rigid $C^*$-tensor category $\text{Deck}^+(S^2|\mathbb{RP}^2)$.

\begin{notation} 
Let $\displaystyle S_{k,n}:=\left\{E \in\mathcal{B}(\CC^n):E^2=E=E^*,
\text{ rank}(E)=k\right\}$ for $0\leq k\leq n$.
\end{notation} 

Note that, as a subspace of $\mathcal{B}(\CC^n)$, $S_{k,n}$ is always compact and in fact homeomorphic to $\text{Gr}_k(n,\CC)$, a complex Grassmannian manifold. Also observe that $S_{0,n}=\{0\}$ and $S_{n,n}=\{I_{n\times n}\}$. 

Our first result yields a formula to express matrix entries of any object $\sigma$ in $\text{Deck}^+(S^2|\mathbb{RP}^2)$.

    \begin{lem}\label{formulalemmarp2}
        Let $\sigma=(\CC^n,u^\sigma)$ be an object. Then, we can write the matrix entries of $u$, indexed by $\tilde{x},\tilde{y}\in S^2$, as \[u_{\tilde{x},\tilde{y}}^\sigma=\delta_{\tilde{y},\tilde{x}}P(\tilde{x})+\delta_{\tilde{y},-\tilde{x}}(I-P(\tilde{x}))\] for a uniquely determined function $P\in C(S^2,S_{k,n})$ satisfying $P(\tilde{x})=P(-\tilde{x})$ for all $\tilde{x}\in S^2$.
    \end{lem}

    \begin{proof} 
        Let $\tilde{x},\tilde{y}\in S^2$. Note that $p(\tilde{x})=p(\tilde{y})$ if and only if $\tilde{y}=\pm\tilde{x}$. Since a quantum deck transformation acts on fibers only, we can write \[u_{\tilde{x},\tilde{y}}^\sigma=\delta_{\tilde{y},\tilde{x}}P(\tilde{x})+\delta_{\tilde{y},-\tilde{x}}Q(\tilde{x})\] for some uniquely determined $P\in C(S^2,S_{k,n})$, $Q\in C(S^2,S_{l,n})$, for some $0\leq k,l\leq n$. 

        Fix $\tilde{x}\in S^2$. Then, since $\sigma$ is a permutation, \[I=\sum_{\tilde{y}\in S^2}u_{\tilde{x},\tilde{y}}^\sigma=\sum_{\tilde{y}\in S^2}\delta_{\tilde{y},\tilde{x}}P(\tilde{x})+\delta_{\tilde{y},-\tilde{x}}Q(\tilde{x})=P(\tilde{x})+Q(\tilde{x}).\] This implies that $Q(\tilde{x})=I-P(\tilde{x})$. 

    Fix $\tilde{y}\in S^2$. Then, similarly, we must have \[I=\sum_{\tilde{x}\in S^2}u_{\tilde{x},\tilde{y}}^\sigma=\sum_{\tilde{x}\in S^2}\delta_{\tilde{y},\tilde{x}}P(\tilde{x})+\delta_{\tilde{y},-\tilde{x}}Q(\tilde{x})=P(\tilde{y})+Q(-\tilde{y}).\] Thus, for all $\tilde{z}\in S^2$, \[P(-\tilde{z})=I-Q(-\tilde{z})=P(\tilde{z}).\] This proves the claim.
    \end{proof}

    \begin{rmk}
        \label{descentrmk}
        Note that $P$ above may be regarded as a continuous function on $\mathbb{RP}^2$ by the condition $P(\tilde{x})=P(-\tilde{x})$ for all $\tilde{x}\in S^2$. Therefore, in the following, we abuse notation denoting $P(\tilde{x})$ by  $P^\sigma(x)$ for 
        $x \in \mathbb{RP}^2$ corresponding to $\sigma$. We also write $\Phi:\sigma\to P^\sigma$.

        We show in \Cref{RP2theorem} below that the classification of $\sigma$ is reduced to the classification of $\Phi(\sigma)$.  
    \end{rmk}

We fix some notation now before stating the first main result. 

\begin{notation}
\begin{enumerate}[{\rm (i)}]

\item For based topological spaces $(X,x_0)$, $(Y,y_0)$, we denote by $C((X,x_0)$, $(Y,y_0))$ the continuous functions which map $x_0$ to $y_0$. We refer to these as ``basepoint-preserving''.

\item If $\Hl$ is a Hilbert space with inner product $(\cdot,\cdot)$ and $\xi,\eta\in\Hl$, let $t_{\xi,\eta}$ be the rank 1 operator defined by $t_{\xi,\eta} \zeta =(\zeta ,\eta)\xi$ for each $\zeta \in\Hl$. Note that $t_{\xi,\eta}^*=t_{\eta,\xi}$. 

\item Let $\{e_i\}_{i=1}^n$ denote the standard orthonormal basis of $\CC^n$. Fix the basepoint $B_k=\sum_{i=1}^kt_{e_i,e_i}\in S_{k,n}$ and let $x_0\in\mathbb{RP}^2$ be a basepoint.  

\item Denote by $U_m$ the unitary group, and define the action of $u\in U_k\oplus U_{n-k}\subseteq U_n$ on $\\P\in C((\mathbb{RP}^2,x_0),(S_{k,n},B_k))$ by $P\mapsto uPu^*$, where $(uPu^*)(x):=uP(x)u^*$ for each $x\in\mathbb{RP}^2$.
\end{enumerate}
\end{notation}

    \begin{rmk}
        \label{classicalRP2}
    For $n=1$, in view of \Cref{formulalemmarp2}, 
    the 1-dimensional classical quantum deck transformations are the usual 
    deck transformations of $S^2|\mathbb{RP}^2$ corresponding to elements of the fundamental group $\pi_1(\mathbb{RP}^2)$. For $n \ge 2$, we have the following: 
    \end{rmk}

\begin{thm}\label{RP2theorem} 
Let $n\geq2$ and $1\leq k<n$. 
We have the following characterization of the category 
$\Deck^+(S^2|\mathbb{RP}^2):$
\begin{enumerate}[{\rm (i)}]
    \item 
    The mapping $\Phi$ 
    above 
    establishes a one-to-one correspondence between 
    unitary equivalence classes of objects $\sigma=(\Hl_\sigma,u^\sigma)$ of $\Deck^+(S^2|\mathbb{RP}^2)$ with $\dim\Hl_\sigma=n$ and $\text{rank}(u_{\tilde{x},\tilde{x}}^\sigma)=k$
    and elements of $C((\mathbb{RP}^2,x_0),(S_{k,n},B_k))/(U_k\oplus U_{n-k})$.

\item Irreducible objects represented on $\CC^n$ correspond to $P\in C((\mathbb{RP}^2,x_0),(S_{k,n},B_k))$ such that the $C^*$-algebra $C^*(P(\mathbb{RP}^2))$ generated by the set $\{P(x):x\in\mathbb{RP}^2\}$ is $M_n(\CC)$.
    \item There is an irreducible object for any such $k$ and $n$.
\end{enumerate}
\end{thm}
\begin{proof}

   For (i), let $\sigma$ and $\tau$ be two objects represented on $\CC^n$ and 
 $P^{\sigma}, P^\tau\in C(\mathbb{RP}^2,S_{k,n})$ be the corresponding functions obtained for $\sigma$ and $\tau$ as in \Cref{descentrmk} respectively. Now, there are unitaries $U,V$ with $UP^\sigma(x_0)U^*=B_k=VP^\tau(x_0)V^*$. By definition of unitary equivalence and \Cref{formulalemmarp2}, the objects $\tau$ and $\sigma$ are unitarily equivalent if and only if there is a unitary operator $T$ with $TP^\sigma(z)T^*=P^\tau(z)$ for all $z\in\mathbb{RP}^2$. 
 
 Thus, $\sigma$ is unitarily equivalent to an object $\sigma'$ with $P^{\sigma'}(x_0)=B_k$, and similarly for $\tau$. So without loss of generality, we can restrict to basepoint-preserving continuous functions $P^\sigma,P^\tau\in C((\mathbb{RP}^2,x_0),(S_{k,n},B_k))$ and unitary operators which fix $B_k$ under conjugation as intertwiners. Then, $B_k(\CC^n)$ is an invariant subspace of $U$. Such a unitary operator decomposes into a direct sum of $U=U_1\oplus U_2\in U_k\oplus U_{n-k}$. This yields (i).

    For (ii), let $\sigma$ be an object, and let $T:\CC^n\rightarrow\CC^n$ be any endomorphism of $\sigma$. By definition, $T$ commutes with all the projections $P^\sigma(x)$, and hence is in the commutant of the $C^*$-algebra generated by 
    $P^\sigma (x) \; (x \in\mathbb{RP}^2)$, denoted by 
    $C^*(P^\sigma (x):x \in\mathbb{RP}^2)$. Similarly, any operator in this commutant $C^*(P^\sigma (x):x\in\mathbb{RP}^2)'$ will be an endomorphism of $\sigma$. Thus, we have $\text{End}(\sigma)=C^*(P^\sigma (x):x\in\mathbb{RP}^2)'$. Then, (ii) follows by definition, as $C^*(P^\sigma (x):x\in\mathbb{RP}^2)$ is a finite-dimensional $C^*$-algebra.

    Finally, for (iii), we have the following construction: 
    
    Fix $k$. Pick any continuous surjection $F:\mathbb{RP}^2\twoheadrightarrow[0,1]$. Using $F$, since $U_n$ is path-connected, we just need to specify a finite collection $\Ul$ of unitary operators such that the $C^*$-algebra $A$ generated by the conjugates $\{U B_kU^*:U\in\Ul\}$ coincides with $M_{n}(\CC)$ (by (ii)). 

    Conjugating $B_k$ by permutations, we obtain all the rank $k$ projections onto subspaces of $\CC^n$ spanned by a subset of $\{e_l:l=1,\dots,n\}$ of size $k<n$. We add these permutations to $\Ul$ so that all these rank $k$ projections are in $A$. Observe that by multiplying pairs of these rank $k$ projections having $k-1$ common rank 1 subprojections, we obtain in $A$ all rank $k-1$ projections to the subspaces of $\CC^n$ spanned by a subset of $\{e_l:l=1,\dots,n\}$ of size $k-1$. Inductively, we have that $A$ contains all the rank $1$ projections $t_{e_l,e_l}$, $l=1,\dots,n$. 

    Now, for fixed $i, j$, let $\xi_{i,j}^{\pm}:=\frac{1}{\sqrt{2}}(e_i\pm e_j)$. We can construct a unitary which maps the basis $\{e_l:l=1,\dots,n\}$ to a basis with $e_i$ and $e_j$ replaced with $\xi_{i,j}^+$ and $\xi_{i,j}^-$, respectively. By adding this unitary to $\Ul$ and repeating the induction above, we can arrange for $A$ to contain all rank 1 projections in this new basis. In particular, $t_{\xi_{i,j}^\pm,\xi_{i,j}^\pm}\in A$. 

    By a routine computation, $\sqrt{2}\cdot t_{\xi_{i,j}^+,\xi_{i,j}^+}t_{e_i,e_i}=t_{\xi_{i,j}^+,e_i}\in A$ for all $i,j=1,\dots,n$. Similarly, $t_{\xi_{i,j}^+,e_j}\in A$. Taking adjoints, both  $t_{e_i,\xi_{i,j}^+}$ and $t_{e_j,\xi_{i,j}^+}\in A$. Then, for all $i$ and $j$ from 1 to $n$, $t_{e_i,e_j}=t_{e_i,\xi_{i,j}^+}\cdot t_{\xi_{i,j}^+,e_j}\in A$. Thus, $A$ contains all matrix units, which generate $M_n(\CC)$. It follows from (ii) that this constructed $\sigma$ is irreducible.

\end{proof}

\section{\texorpdfstring{$\pi_1^q$}{} of General Spaces with Finite-Sheeted Coverings} \label{sec:qfgfscsp}

With a basic case covered, we prove a generalization of \Cref{RP2theorem}. We omit detailed proofs in favor of highlighting which modifications to the statements and arguments above must be made. We have the following setting:
    
    Let $X$ be a path-connected, locally compact Hausdorff space with $\pi_1(X)=:\Gamma$ of order $N>1$. Moreover, assume $X$ is locally path-connected and semi-locally simply connected. Let $\widetilde{X}$ be the universal cover with covering map $p:\widetilde{X}\rightarrow X$, and let $\alpha:\Gamma\rightarrow\text{Deck}(\widetilde{X})$ denote the 
    left action of $\Gamma$ on $\widetilde{X}$ by (classical) deck transformations. 
    Then, $X\cong\widetilde{X}/\Gamma$. We denote by $e$ the neutral element of $\Gamma$. For convenience later, we will also fix an enumeration of the nonidentity elements of $\Gamma$ as $\{g_1,g_2,\dots,g_{N-1}\}=\Gamma-\{e\}$.

    \begin{lemma}\label{formulalemmagen}
        Let $\mathcal{H}$ be a finite-dimensional Hilbert space. Objects $\sigma=(\mathcal{H},u)$ correspond by the formula \[u_{\tilde{x},\tilde{y}}^\sigma=\sum_{g\in\Gamma}\delta_{\tilde{y},\alpha(g)(\tilde{x})}P_g(\tilde{x})\] uniquely to tuples $(P_g)_{g\in\Gamma-\{e\}}$, with each $P_g\in C(\widetilde{X},S_{k_g,n})$ for some $0\leq k_g\leq n$ and $P_e=I-\sum_{g\in\Gamma-\{e\}}P_g$ such that, if $g,h\in\Gamma$ are distinct, $P_g(\tilde{x})P_h(\tilde{x})=0$ for all $\tilde{x}\in\widetilde{X}$, and 
        \begin{equation}\label{equivariance}
        \sum_{g\in\Gamma-\{e\}}P_g(\tilde{x})=\sum_{g\in\Gamma-\{e\}}P_g(\alpha(g^{-1})(\tilde{x}))
        \end{equation} 
        for $\tilde{x}\in\widetilde{X}$. 
    
    \end{lemma}
\begin{rmk}\label{matrixrmk}
    The condition \Cref{equivariance} above is equivalent to the condition that the matrix \[\begin{bmatrix}
        P_{g^{-1}h}(\alpha(g)(\tilde{x}))
    \end{bmatrix}_{(g,h)\in\Gamma\times\Gamma}\] define a quantum permutation on a set of size $N$ for each $\tilde{x}$. That is, the functions $P_g$ must define a family of representations of $C(S_{N}^+)$ as above, indexed by points in the cover $\widetilde{X}$.
\end{rmk}
    \begin{proof}
Let $\sigma=(\mathcal{H},u^\sigma)$ be an object. The proof in \Cref{formulalemmarp2} carries over almost exactly, but now we use that $p(\tilde{x})=p(\tilde{y})$ if and only if there exists a $g\in\Gamma$ such that $\tilde{y}=\alpha(g)(\tilde{x})$. Since $u^\sigma$ acts on the fibers only, we write 
\[u_{\tilde{x},\tilde{y}}^\sigma=\sum_{g\in\Gamma}\delta_{\tilde{y},\alpha(g)(\tilde{x})}u_{\tilde{x},\alpha(g)(\tilde{x})}^\sigma.
\]  

        We obtain the functions 
        $P_g(\tilde{x}):=u_{\tilde{x},\alpha(g)(\tilde{x})}^\sigma$. 
        We have 
\[I=\sum_{\tilde{y}\in\widetilde{X}}u_{\tilde{x},\tilde{y}}^\sigma=\sum_{g\in\Gamma}P_g(\tilde{x}).
        \]
        This gives pairwise orthogonality for the functions $P_g$ pointwise.
        
        Also, for any $\tilde{y}\in\widetilde{X}$, 
        \begin{align*}
            I&=\sum_{\tilde{x}\in\widetilde{X}}u_{\tilde{x},\tilde{y}}^\sigma
            =\sum_{\tilde{x}\in\widetilde{X}}\sum_{g\in\Gamma}\delta_{\tilde{y},\alpha(g)(\tilde{x})}P_g(\tilde{x})\\
            &=\sum_{\tilde{x}\in\widetilde{X}}\sum_{g\in\Gamma}\delta_{\alpha(g^{-1})(\tilde{y}),\tilde{x}}P_g(\alpha(g^{-1})(\tilde{y})) 
            =\sum_{g\in\Gamma}P_g(\alpha(g^{-1})(\tilde{y})). 
        \end{align*}

        Then, setting these equal for a chosen $\tilde{x}=\tilde{y}=\tilde{z}$ and subtracting the function $P_e$ from both sides, we obtain the relation \[\sum_{g\in\Gamma-\{e\}}P_g(\tilde{z})=\sum_{g\in\Gamma-\{e\}}P_g(\alpha(g^{-1})(\tilde{z})).\]

        Note that, if we are just given projections for non-identity elements of $\Gamma$, we reconstruct $P_e=I-\sum_{g\neq e}P_g$.
    \end{proof}
    As above, we use superscripts to label, for an object $\sigma$, the associated tuple of functions $(P_g^\sigma)_{g\in\Gamma-\{e\}}$.

    Now, we characterize $\text{Deck}^+(\widetilde{X}|X)$ as we did for $X=\mathbb{RP}^2$ and $\widetilde{X}=S^2$ in \Cref{RP2theorem}. We first introduce some notation:

    Let $\NN$ denote the non-negative integers. For $n\geq2$, define the set 
    \[
    L_{N,n}:=\bigg\{(k_g)_{g\in\Gamma-\{e\}}\in\NN^{\Gamma-\{e\}}:0<\sum_{g\in\Gamma-\{e\}}k_g\leq n
    \bigg\}.
    \]
    
    Note that $(k_g)_{g\in\Gamma-\{e\}}\in L_{N,n}$ always contains a nonzero component.

    For $\vec{k}\in L_{N,n}$, let $S_{\vec{k},n}$ denote the set \[\{(P_{g})_{g\in{\Gamma-\{e\}}}\in\mathcal{B}(\CC^n)^{\Gamma-\{e\}}:P_g^2=P_g^*=P_g,\text{rk}(P_g)=k_g,P_gP_h=\delta_{g,h}P_g, \,  g,h\in\Gamma-\{e\}\}.\] This is the space of $(N-1)$-tuples of projections in $\mathcal{B}(\CC^n)$ which are pairwise mutually orthogonal and such that the $i^{\text{th}}$ projection has prescribed rank $k_i$. 
    
    This may also be interpreted as a flag manifold, where we identify the increasing sequence of projections $P_{g_1}\subseteq P_{g_1}+P_{g_2}\subseteq\dots\subseteq P_{g_1}+P_{g_2}+\dots+P_{g_{N-1}}$ with the flag given by the associated increasing sequence of subspaces of $\CC^n$. In particular, if $N=2$, we obtain a Grassmannian manifold. The signature is determined from the vector $\vec{k}$ as $(k_{g_1},k_{g_1}+k_{g_2},\dots,k_{g_1}+\dots+k_{g_{N-1}})$.

    Let $n\geq2$ and $\vec{k}\in L_{N,n}$. We pick a basepoint \[B_{\vec{k}}:=\bigg(\sum_{i=1}^{k_{g_1}}t_{e_i,e_i},\sum_{i=1}^{k_{g_2}}t_{e_{i+k_{g_1}},e_{i+k_{g_1}}},\dots,\sum_{i=1}^{k_{g_{N-1}}}t_{e_{i+k_{g_1}+\dots+k_{g_{N-2}}},e_{i+k_{g_1}+\dots+k_{g_{N-2}}}}\bigg)\in S_{\vec{k},n}.\] That is, $B_{\vec{k}}$ is a tuple of projections with the first entry the rank $k_{g_1}$ projection to the first $k_{g_1}$ basis vectors, the second entry the rank $k_{g_2}$ projection to the next $k_{g_2}$ basis vectors, and so on. We also fix a basepoint $\tilde{x}_0\in\widetilde{X}$.

    \begin{rmk}
        \label{classicalGen}
    For $n=1$, in view of \Cref{formulalemmagen}, 
    the 1-dimensional classical quantum deck transformations are the usual 
    deck transformations of $\widetilde{X}|X$ corresponding to elements of the fundamental group $\pi_1(X)$. For $n \ge 2$, we have the following: 
    \end{rmk}

    \begin{thm}\label{generaltheorem}
        We have the following description of the category of quantum deck transformations associated to the universal cover $\widetilde{X}$: 
        \begin{enumerate}[{\rm (i)}]
       
        \item Let $\sigma=(\CC^n,u)$ be an object. The unitary equivalence class of $\sigma$ corresponds with an element in the quotient \[\bigg\{P=(P_g)_{g\in\Gamma-\{e\}}\in C((\widetilde{X},\tilde{x}_0),(S_{\vec{k},n},B_{\vec{k}})):\sum_{g\in\Gamma-\{e\}}P_g=\sum_{g\in\Gamma-\{e\}}P_g\circ\alpha(g^{-1})\bigg\}/G\] for some $\vec{k}\in L_{N,n}$ under the action of \[G:=U_{k_{g_1}}\oplus U_{k_{g_2}}\oplus\cdots\oplus U_{k_{g_{N-1}}}\oplus U_{n-k_{g_1}-\dots-k_{g_{N-1}}}\hookrightarrow U_n\] given by $(P_g)_{g\in\Gamma-\{e\}}\in S_{\vec{k},n}\mapsto (uP_gu^*)_{g\in\Gamma-\{e\}}$ for $u\in G$.

        \item An object $\sigma$ with associated function $P\in C((\widetilde{X},\tilde{x}_0),(S_{\vec{k},n},B_{\vec{k}}))$ (as in  \Cref{formulalemmagen}) is irreducible if and only if the $C^*$-algebra generated by the set $\{P_g(\widetilde{X})\}_{g\in\Gamma-\{e\}}\subseteq \mathcal{B}(\mathcal{H})$, denoted $C^*(P_{g_1}(\tilde{z}),\dots,P_{g_{N-1}}(\tilde{z}):\tilde{z}\in\widetilde{X})$, equals $\mathcal{B}(\CC^n)\cong M_n(\CC)$. 
        \item There exist irreducible (quantum) objects for any $n\geq2$ and $\vec{k}\in L_{N,n}$ with all components of $\vec{k}$ strictly less than $n$.
        \end{enumerate}
    \end{thm}
\begin{proof}

For (i) and (ii), we can essentially repeat the arguments for the analogous statements in \Cref{RP2theorem} verbatim. Similarly, (iii) can be proven by applying the process in the proof of (iii) in \Cref{RP2theorem} to a nonzero projection in the tuple $(P_g)_{g\in\Gamma-\{e\}}$. The only obstruction here is that we must be careful to ensure that our surjection $F:\widetilde{X}\twoheadrightarrow[0,1]$ is compatible with the equivariance condition \eqref{equivariance} on the projections obtained along the path. This can be done as follows:

    Let $n\geq2$, $\vec{k}\in L_{N,n}$. We have the following sequence of surjective, continuous maps: \[\widetilde{X}\twoheadrightarrow \widetilde{X}/\Gamma \cong X\twoheadrightarrow[0,1],\] where the second map is obtained by Urysohn's lemma and connectedness of $X$. 
    
    Let $F:\widetilde{X}\twoheadrightarrow[0,1]$ denote the composition of these maps. Then, note that $F$ satisfies $F=F\circ\alpha(g)$ for all $g\in\Gamma$. This now permits the reduction to finding a continuous path $[0,1]\rightarrow S_{\vec{k},n}$, as condition \eqref{equivariance} is automatically satisfied: given a path $\hat{P_g}:[0,1]\rightarrow S_{k_g,n}$, $g\neq e$, we lift to $P_g:\widetilde{X}\rightarrow S_{k_g,n}$ by $P_g(\tilde{x})=\hat{P_g}(F(\tilde{x}))$. Then, for any $g,h\in\Gamma$, \[P_g(\alpha(h)(\tilde{x}))=\hat{P_g}(F(\alpha(h)(\tilde{x})))=\hat{P_g}(F(\tilde{x}))=P_g(\tilde{x}).\]
    
\end{proof}

We conclude this section with a short remark:

\begin{rmk} \label{representations}
    We can adopt the following perspective on a quantum deck transformations of a space with an $N$-sheeted universal cover: a quantum deck transformation is a continuous family of representations $\rho_{\tilde{x}}:C(S_{N}^+)\rightarrow \Bl(\CC^n)$ 
    indexed by $\tilde{x}\in \widetilde{X}$, the universal cover of a space $X$. 
    Here, $S_N^+$ are the compact quantum permutation groups on $N$ elements of \cite{Wan98}, and $C(S_N^+)$ is the associated algebra of functions on $S_N^+$.
    
    Consider the case where the cover $\widetilde{X}|X$ is two-sheeted. It is interesting to note that we cannot have irreducibility of the representations $\rho_{\tilde{x}}:C(S_2^+)\rightarrow \mathcal{B}(\CC^n)$ at each point 
    $\tilde x\in \widetilde{X}$. This is a consequence of elementary Gelfand-Naimark theory applied to the commutative $C^*$-algebra $C(S_2^+)$ together with the connectedness of $\widetilde{X}$. A similar argument would work for any space with a 3-sheeted universal cover, since $C(S_3^+)$ is also commutative. 
    
    It would be interesting to determine whether irreducible objects exist which arise from pointwise-irreducible families of representations in $N$-sheeted cases, $N\geq4$. We consider this perspective further in the next section.
\end{rmk}

\section{Classicality}\label{sec:clas}

Having constructed the examples above, we ask the following two natural questions. First, are these discrete quantum groups genuinely quantum (i.e. their function algebras are both noncommutative and non-cocommutative)? Second, do any properties of the classical fundamental group hold in this quantum setting? In this section, we investigate the first question, and the next section, the second question. 

\Cref{RP2theorem} characterizes $\pi_1^q(\mathbb{RP}^2)$. Despite the natural construction as a discrete quantum symmetry group with nonclassical deck transformations, this example is classical in the following sense.

\begin{thm}\label{thm:qfpz2ccm}
	$\pi_1^q(\mathbb{RP}^2)$ is the dual discrete quantum group to a compact group.
\end{thm}

Before the proof of this result, we fix some notation. 

\begin{notation}
\begin{enumerate}[\rm(i)]
\item 
For $C^*$-algebras $A$ and $B$, define $\phi: A \otimes B \to B \otimes A$ by
 $\phi(a\otimes b)=b\otimes a$ for all $a \in A, b\in B$.   
 
\item 
For Hilbert spaces $\Hl$ and $\Kl$,
	define the ``flip'' $t:\Hl \otimes\Kl \to \Kl \otimes\Hl$ by 
    $t(\xi\otimes\eta)=\eta\otimes\xi$, $\xi \in \Hl, \eta \in \Kl$. 
\end{enumerate}
\end{notation}

Note that 
$t^*$ is another ``flip'', $t^*(\eta\otimes\xi) = \xi\otimes\eta$ in the notation above, and that $t^{-1}=t^*$, so $t$ is unitary. 
    If $\Hl = \Kl$, then $t^{-1}=t$. Furthermore, if $A \subset \Bl(\Hl)$ and  
     $B \subset \Bl(\Kl)$ are $C^*$-subalgebras, 
     then $\phi(a\otimes b) = t (a \otimes b) t^*$.
     
\begin{proof}
	It suffices to show $C_0(\pi_1^q(\mathbb{RP}^2))$ is cocommutative, i.e., $\phi\circ\Delta=\Delta$, where 
    $\phi$ is as defined above with $A=B=C_0(\pi_1^q(\mathbb{RP}^2))$. Then, the product given by $\omega_1*\omega_2:=(\omega_1\otimes\omega_2)\circ\Delta$ for $\omega_1,\omega_2$ finitely-supported linear functionals on $C_0(\pi_1^q(\mathbb{RP}^2))$, makes the dual multiplier Hopf *-algebra \cite{Vdl-duality98} into a commutative algebra. 
	
	Let $\sigma=(\Hl_\sigma,u^\sigma)$ and $\tau=(\Hl_\tau,u^\tau)$ be objects in $\Deck^+(S^2|\mathbb{RP}^2)$. 
	By \Cref{formulalemmarp2} and \Cref{RP2theorem}, we can write \[u^\sigma_{\tilde{x},\tilde{y}}=\delta_{\tilde{y},\tilde{x}}P_1(\tilde{x})+\delta_{\tilde{y},-\tilde{x}}P_2(\tilde{x})\] and \[u^\tau_{\tilde{x},\tilde{y}}=\delta_{\tilde{y},\tilde{x}}Q_1(\tilde{x})+\delta_{\tilde{y},-\tilde{x}}Q_2(\tilde{x}),\] where $P_2(\tilde{x})=I_{\Hl_\sigma}-P_1(\tilde{x})$ and $Q_2(\tilde{x})=I_{\Hl_\tau}-Q_1(\tilde{x})$, and $P_1(-\tilde{x})=P_1(\tilde{x})$ and $Q_1(-\tilde{x})=Q_1(\tilde{x})$, for all $\tilde{x}\in S^2$.
	
	Let $\tilde{x},\tilde{y}\in S^2$. We have that 
	\begin{align*}
		u^{\sigma\boxtimes\tau}_{\tilde{x},\tilde{y}}
		&= \sum_{\tilde{z}\in S^2}u^\sigma_{\tilde{x},\tilde{z}}\otimes u^\tau_{\tilde{z},\tilde{y}}\\
		&=\delta_{\tilde{y},\tilde{x}}(P_1(\tilde{x})\otimes Q_1(\tilde{y})+P_2(\tilde{x})\otimes Q_2(-\tilde{y}))+\delta_{\tilde{y},-\tilde{x}}(P_2(\tilde{x})\otimes Q_1(\tilde{y})+P_1(\tilde{x})\otimes Q_2(-\tilde{y}))\\
		&=\delta_{\tilde{y},\tilde{x}}(P_1(\tilde{x})\otimes Q_1(\tilde{y})+P_2(\tilde{x})\otimes Q_2(\tilde{y}))+\delta_{\tilde{y},-\tilde{x}}(P_2(\tilde{x})\otimes Q_1(\tilde{y})+P_1(\tilde{x})\otimes Q_2(\tilde{y}))
		\end{align*} 
		where $Q_2(-\tilde{y})=Q_2(\tilde{y})$ is used, and 
		\begin{align*}
		u^{\tau\boxtimes\sigma}_{\tilde{x},\tilde{y}}
		&=\sum_{\tilde{z}\in S^2}u^\tau_{\tilde{x},\tilde{z}}\otimes u^\sigma_{\tilde{z},\tilde{y}}\\
		&=\delta_{\tilde{y},\tilde{x}}(Q_1(\tilde{x})\otimes P_1(\tilde{y})+Q_2(\tilde{x})\otimes P_2(-\tilde{y}))+\delta_{\tilde{y},-\tilde{x}}(Q_2(\tilde{x})\otimes P_1(\tilde{y})+Q_1(\tilde{x})\otimes P_2(-\tilde{y}))\\
		&=\delta_{\tilde{y},\tilde{x}}(Q_1(\tilde{x})\otimes P_1(\tilde{y})+Q_2(\tilde{x})\otimes P_2(\tilde{y}))+\delta_{\tilde{y},-\tilde{x}}(Q_2(\tilde{x})\otimes P_1(\tilde{y})+Q_1(\tilde{x})\otimes P_2(\tilde{y}))
	\end{align*}
	where $P_2(-\tilde{y})=P_2(\tilde{y})$ is used.
	
	Let $t:\Hl_\sigma\otimes\Hl_\tau\rightarrow\Hl_\tau\otimes\Hl_\sigma$ be the ``flip'' defined above. 
	We have that \begin{align*}
		tu^{\sigma\boxtimes\tau}_{\tilde{x},\tilde{y}}t^*
		&=\delta_{\tilde{y},\tilde{x}}(Q_1(\tilde{y})\otimes P_1(\tilde{x})+Q_2(\tilde{y})\otimes P_2(\tilde{x}))
		+\delta_{\tilde{y},-\tilde{x}}(Q_1(\tilde{y})\otimes P_2(\tilde{x})+Q_2(\tilde{y})\otimes P_1(\tilde{x}))
		\\
		&=\delta_{\tilde{y},\tilde{x}}(Q_1(\tilde{x})\otimes P_1(\tilde{y})+Q_2(\tilde{x})\otimes P_2(\tilde{y}))
		+\delta_{\tilde{y},-\tilde{x}}(Q_1(-\tilde{x})\otimes P_2(-\tilde{y})+Q_2(-\tilde{x})\otimes P_1(-\tilde{y}))
		\\
		&=\delta_{\tilde{y},\tilde{x}}(Q_1(\tilde{x})\otimes P_1(\tilde{y})+Q_2(\tilde{x})\otimes P_2(\tilde{y}))+\delta_{\tilde{y},-\tilde{x}}(Q_1(\tilde{x})\otimes P_2(\tilde{y})+Q_2(\tilde{x})\otimes P_1(\tilde{y}))
		\\
		&=u^{\tau\boxtimes\sigma}_{\tilde{x},\tilde{y}},
		\end{align*}
	where $Q_1(-\tilde{x})=Q_1(\tilde{x})$, 
	$P_2(-\tilde{x})=P_2(\tilde{x})$,  
	$Q_2(-\tilde{x})=Q_2(\tilde{x})$, and
	$P_1(-\tilde{x})=P_1(\tilde{x})$ are used.

	Now, if 
	$\displaystyle u_{\tilde{x},\tilde{y}}=\sum_{\sigma\in\text{Irred}_1^q}u^\sigma_{\tilde{x},\tilde{y}}$ are the universal matrix entries, we have 
	\begin{align*}
		\Delta(u_{\tilde{x},\tilde{y}})
		&=\sum_{\tilde{z}\in S^2}u_{\tilde{x},\tilde{z}}\otimes u_{\tilde{z},\tilde{y}}
        =\sum_{\tilde{z}\in S^2}\sum_{\sigma,\tau\in\text{Irred}_1^q}u^\sigma_{\tilde{x},\tilde{z}}\otimes u^\tau_{\tilde{z},\tilde{y}}\\
		&=\sum_{\sigma,\tau\in\text{Irred}_1^q}u^{\sigma\boxtimes\tau}_{\tilde{x},\tilde{y}}
        =\sum_{\sigma,\tau\in\text{Irred}_1^q}tu^{\tau\boxtimes\sigma}_{\tilde{x},\tilde{y}}t^* 
        \\
        &= \sum_{\sigma,\tau\in\text{Irred}_1^q}\phi(u^{\tau\boxtimes\sigma}_{\tilde{x},\tilde{y}})
        =\phi( \Delta(u_{\tilde{x},\tilde{y}})).
		\end{align*}

    Cocommutativity follows.

\end{proof}

\begin{rmk}
    Note that a similar argument works for any space $X$ with $\pi_1(X)\cong\ZZ/2\ZZ$. Hence, any space with $\pi_1(X)\cong\ZZ/2\ZZ$ has a cocommutative algebra $C_0(\pi_1^q(X))$.
\end{rmk}

Given the above, is $\pi_1^q(X)$ genuinely quantum for any $X$? That is, is there any space  $X$ where $\pi_1^q(X)$ is not dual to a compact group? 
The next theorem gives an affirmative answer. 

First, note that if $C_0(\pi_1^q(X))$ is cocommutative, then for any irreducible  $\sigma,\tau$, and for $\tilde{x},\tilde{y}\in\widetilde{X}$, the universal cover of $X$, \begin{align*}
u^{\sigma\boxtimes\tau}_{\tilde{x},\tilde{y}}
&=\sum_{\tilde{z}\in\widetilde{X}}u^\sigma_{\tilde{x},\tilde{z}}\otimes u^\tau_{\tilde{z},\tilde{y}}
=(e^\sigma\otimes e^\tau)\sum_{\tilde{z}\in\widetilde{X}}u_{\tilde{x},\tilde{z}}\otimes u_{\tilde{z},\tilde{y}}
\\&=(e^\sigma\otimes e^\tau)\Delta(u_{\tilde{x},\tilde{y}})
=(e^\sigma\otimes e^\tau)(\phi\circ\Delta)(u_{\tilde{x},\tilde{y}})
\\&=\phi\bigg((e^\tau\otimes e^\sigma)\sum_{\tilde{z}\in\widetilde{X}}u_{\tilde{x},\tilde{z}}\otimes u_{\tilde{z},\tilde{y}}\bigg)
=\phi(u^{\tau\boxtimes\sigma}_{\tilde{x},\tilde{y}}).
\end{align*}

Therefore, to show $C_0(\pi_1^q(X))$ is non-cocommutative, it suffices to exhibit irreducible quantum deck transformations $\sigma$ and $\tau$ in the associated rigid $C^*$-tensor category such that the flip $t$ does not implement a unitary equivalence $\sigma\boxtimes\tau\sim\tau\boxtimes \sigma$.

\begin{thm}\label{nonclassical}
    $C_0(\pi_1^q(\mathbb{RP}^2\times\mathbb{RP}^2))$ is non-cocommutative.
\end{thm}
    Before we begin the proof, we fix some notation and explain the idea of the argument. 
    
    The universal cover of $\mathbb{RP}^2\times\mathbb{RP}^2$ is $S^2\times S^2$ with product covering map $p\times p$. The group $\pi_1(\mathbb{RP}^2\times\mathbb{RP}^2)\cong\ZZ/2\ZZ\oplus\ZZ/2\ZZ$ acts on the universal cover by 
    \begin{align*}
        \alpha_{(0,1)}(z_1,z_2)&=(z_1,-z_2)\\
        \alpha_{(1,0)}(z_1,z_2)&=(-z_1,z_2)\\
        \alpha_{(1,1)}(z_1,z_2)&=(-z_1,-z_2)
    \end{align*} 
    for $z_1,z_2\in S^2$. 
        For simplicity, denote ${ij}:={(i,j)}$ ($i,j\in\ZZ/2\ZZ$). 

    We need to exhibit a pair of irreducible quantum deck transformations $\sigma$ 
    and $\tau$ for which the flip does not implement a unitary equivalence $\sigma\boxtimes\tau\sim\tau\boxtimes\sigma$. 
   We construct $P_{ij} \in C(S^2 \times S^2,S_{k_{ij},4})$ such that  
$\sum_{i,j \in \ZZ/2\ZZ} P_{ij}(\tilde{x}) = 1$ 
 and 
    $P_{ij}(\alpha_{ij}(\tilde{x}))=P_{ij}(\tilde{x})$ for all $\tilde{x}\in S^2\times S^2$.  
    Then condition \Cref{equivariance} is automatically satisfied, since every element of the Klein group is its own inverse. We obtain a quantum deck transformation $\sigma=(\Hl_\sigma,u^\sigma)$ by the formula \[u^\sigma_{\tilde{x},\tilde{y}}=\delta_{\tilde{y},\tilde{x}}P_{00}(\tilde{x})+\delta_{\tilde{y},\alpha_{10}(\tilde{x})}P_{10}(\tilde{x})+\delta_{\tilde{y},\alpha_{01}(\tilde{x})}P_{01}(\tilde{x})+\delta_{\tilde{y},\alpha_{11}(\tilde{x})}P_{11}(\tilde{x}),\] for all $\tilde{x},\tilde{y}\in S^2\times S^2$, as in \Cref{formulalemmagen}.
    
    On the other hand, we do not require that $P_{ij}$ be invariant under any other classical deck transformations. We combine $\sigma$ with a particular classical deck transformation $\tau$ to obtain the desired non-commutativity of the tensor structure.

\begin{proof}

    Define $P_{01}$ and $P_{11}$ to be 0 everywhere. Pick distinct $x_0,x_1,x_2\in S^2$ which are not antipodal ($p(x_i)\neq p(x_j)$ for any $i,j\in\{0,1,2\}$), and distinct $y_0,y_1,y_2,y_3\in S^2$ such that 
    \[
    \{y_0,y_1,y_2,y_3\}\cap\{-y_0,-y_1,-y_2,-y_3\}=\emptyset.
    \]
    Define a rank-one projection-valued function $P:S^2\rightarrow \Bl(\CC^2)$ such that \begin{enumerate}[\rm (i)]
        \item $P(-z)=P(z)$ for all $z\in S^2$;
        \item $P(x_0)=e_{11}$, the rank-one projection to the subspace $\CC e_1$;  
        \item $P(x_1)=e_{22}$, the rank-one projection to the subspace $\CC e_2$; and
        \item $P(x_2)=t_{\frac{1}{\sqrt{2}}(e_1+e_2),\frac{1}{\sqrt{2}}(e_1+e_2)}$, the rank-one projection to the subspace $\CC\cdot(e_1+e_2)$.
    \end{enumerate}
    We can do this by first defining a function $f:\mathbb{RP}^2\rightarrow[0,1]$ such that $f(p(x_0))=0$, $f(p(x_1))=1/2$ and $f(p(x_2))=1$. Then, we can take a continuous path of unitaries $\{U(t)\}_{t\in[0,1]}\subset\Bl(\CC^2)$ such that $U(0)e_{11}U(0)^*=e_{11}$, $U(1/2)e_{11}U(1/2)^*=e_{22}$, and $U(1)e_{11}U(1)^*=t_{\frac{1}{\sqrt{2}}(e_1+e_2),\frac{1}{\sqrt{2}}(e_1+e_2)}$ by path-connectedness of $U_2$. Then, $P(x):=(U\circ f\circ p)(x)e_{11}(U\circ f\circ p)(x)^*$ satisfies the above.
    
    Next, define a rank-one-projection-valued function $Q:S^2\rightarrow\Bl(\CC^2)$ satisfying the following:\begin{enumerate}[\rm (i)]
        \item $Q(y_0)=e_{11}=Q(-y_0)$,
        \item $Q(y_1)=e_{22}=Q(-y_1)$, 
        \item $Q(y_2)=t_{\frac{1}{\sqrt{2}}(e_1+e_2),\frac{1}{\sqrt{2}}(e_1+e_2)}=Q(-y_2)$,
        \item $Q(y_3)=e_{11}$, and
        \item $Q(-y_3)=e_{22}$.
    \end{enumerate}
    This can be done similarly: take a function on $S^2$ which takes on distinct values $t_1=0<t_2<t_3<\dots<t_8=1$ in $[0,1]$ at each of the eight points in $\{\pm y_0,\pm y_1,\pm y_2,\pm y_3\}$. We can construct a path of unitaries acting on $e_{11}$ by conjugation by concatenating paths defined on the subintervals $[t_i,t_{i+1}]$ in $[0,1]$, each interpolating the desired projections under conjugation. We obtain, by composing, the function $Q:S^2\rightarrow[0,1]\rightarrow\Bl(\CC^2)$.

    Define $P_{10}(z_1,z_2):=P(z_1)\otimes Q(z_2)$ for each $\tilde{x}=(z_1,z_2)\in S^2\times S^2$. Note that $P_{10}\circ\alpha_{10}=P_{10}$. This means the condition \[P_{01}+P_{10}+P_{11}=P_{01}\circ\alpha_{01}+P_{10}\circ\alpha_{10}+P_{11}\circ\alpha_{11}\] is automatic, as $P_{01}$ and $P_{11}$ are just 0. If we set $P_{00}(\tilde{x}):=I-P_{10}(\tilde{x})$ for each $\tilde{x}\in\widetilde{X}$, we obtain a quantum deck transformation $\sigma$, given by the formula \[u^\sigma_{\tilde{x},\tilde{y}}=\delta_{\tilde{y},\tilde{x}}P_{00}(\tilde{x})+\delta_{\tilde{y},\alpha_{10}(\tilde{x})}P_{10}(\tilde{x}).\] Since the range of $P_{10}$ generates $\Bl(\CC^2\otimes\CC^2)\cong M_4(\CC)$ by construction, $\sigma$ is irreducible by \Cref{RP2theorem}.

    Let $\tau=(\CC,u^\tau)$ be the (irreducible) classical deck transformation given by \[u^\tau_{\tilde{x},\tilde{y}}=\delta_{\tilde{y},\alpha_{01}(\tilde{x})}\cdot1.\] Observe that \begin{align*}
    u^{\sigma\boxtimes\tau}_{\tilde{x},\tilde{y}}&=\sum_{\tilde{z}\in S^2\times S^2}u^\sigma_{\tilde{x},\tilde{z}}\otimes u^\tau_{\tilde{z},\tilde{y}}=\sum_{\tilde{z}\in S^2\times S^2}\delta_{\tilde{z},\tilde{x}}\delta_{\tilde{y},\alpha_{01}(\tilde{z})}P_{00}(\tilde{x})\otimes 1+\delta_{\tilde{z},\alpha_{10}(\tilde{x})}\delta_{\tilde{y},\alpha_{01}(\tilde{z})}P_{10}(\tilde{x})\otimes 1\\
    &=\,\delta_{\tilde{y},\alpha_{01}(\tilde{x})}P_{00}(\tilde{x})\otimes 1+\delta_{\tilde{y},\alpha_{11}(\tilde{x})}P_{10}(\tilde{x})\otimes 1,\end{align*} and \begin{align*}
        u^{\tau\boxtimes\sigma}_{\tilde{x},\tilde{y}}&=\sum_{\tilde{z}\in S^2\times S^2}u^\tau_{\tilde{x},\tilde{z}}\otimes u^\sigma_{\tilde{z},\tilde{y}}=\sum_{\tilde{z}\in S^2\times S^2}\delta_{\tilde{z},\alpha_{01}(\tilde{x})}\delta_{\tilde{y},\tilde{z}}1\otimes P_{00}(\tilde{z})+\delta_{\tilde{z},\alpha_{01}(\tilde{x})}\delta_{\tilde{y},\alpha_{10}(\tilde{z})}1\otimes P_{10}(\tilde{z})\\
       & =\,\delta_{\tilde{y},\alpha_{01}(\tilde{x})}1\otimes P_{00}(\alpha_{01}(\tilde{x}))+\delta_{\tilde{y},\alpha_{11}(\tilde{x})}1\otimes P_{10}(\alpha_{01}(\tilde{x})).
    \end{align*}
    
    Suppose $\sigma\boxtimes\tau$ and $\tau\boxtimes\sigma$ are unitarily equivalent. By the usual identifications $\Hl_\sigma\otimes\CC\cong\Hl_\sigma\cong\CC\otimes\Hl_\sigma$, we see that $\sigma\boxtimes\tau$ is unitarily equivalent to $\sigma'=(\CC^2\otimes\CC^2,u^{\sigma'})$, with \[u^{\sigma'}_{\tilde{x},\tilde{y}}:=\delta_{\tilde{y},\alpha_{01}(\tilde{x})}P_{00}(\tilde{x})+\delta_{\tilde{y},\alpha_{11}(\tilde{x})}P_{10}(\tilde{x}).\] Similarly, $\tau\boxtimes\sigma$ is unitarily equivalent to $\sigma''=(\CC^2\otimes\CC^2,u^{\sigma''})$ with 
    \[u^{\sigma''}_{\tilde{x},\tilde{y}}:=\delta_{\tilde{y},\alpha_{01}(\tilde{x})}P_{00}(\alpha_{01}(\tilde{x}))+\delta_{\tilde{y},\alpha_{11}(\tilde{x})}P_{10}(\alpha_{01}(\tilde{x})).\]
    By assumption, ${\sigma'}$ and $\sigma''$ are unitarily equivalent. That is, there exists a unitary $U:\CC^2\otimes\CC^2\rightarrow\CC^2\otimes\CC^2$ such that, for all $\tilde{x},\tilde{y}\in S^2\times S^2$, \[Uu^{\sigma''}_{\tilde{x},\tilde{y}}U^*=u^{{\sigma'}}_{\tilde{x},\tilde{y}}.\] In particular, for any $\tilde{x}=(z_1,z_2)$, by setting $\tilde{y}=\alpha_{11}(\tilde{x})$, we must have \[UP_{10}(z_1,-z_2)U^*=P_{10}(z_1,z_2)\] by a slight abuse of notation. By construction of $P_{10}$, choosing $z_1$ and $z_2$ appropriately, this implies that $U$ fixes the rank one projections $e_{11}\otimes e_{11}$, $e_{11}\otimes e_{22}$, $e_{22}\otimes e_{11}$, and $e_{22}\otimes e_{22}$ under conjugation. 
    
    Thus, each of the subspaces $\CC e_1\otimes e_1$, $\CC e_1\otimes e_2$, $\CC e_2\otimes e_1$, and $\CC e_2\otimes e_2$ are invariant under $U$, so $U$ is a diagonal operator. But we must also have \[U(e_{11}\otimes e_{22})U^*=UP(x_0,-y_2)U^*=P(x_0,y_2)=e_{11}\otimes e_{11}.\] Since $U$ is diagonal, this implies $e_{11}=e_{22}$, a contradiction. 

    This shows that $\sigma\boxtimes\tau$ and $\tau\boxtimes\sigma$ cannot be unitarily equivalent. In particular, the flip $t$ cannot implement a unitary equivalence, and the result follows.
\end{proof}

The above theorem implies a result of some theoretical importance. Recall that for discrete quantum groups $\Gamma_1$ and $\Gamma_2$ with comultiplications $\Delta_1$ and $\Delta_2$, respectively, their product is $(B,\Delta)$, where \[B=C_0(\Gamma_1\times\Gamma_2):=C_0(\Gamma_1)\otimes C_0(\Gamma_2)\] and \[\Delta:=\sigma_{23}\circ(\Delta_1\otimes\Delta_2).\] Here $\sigma_{23}(a_1\otimes a_2\otimes b_1\otimes b_2)=a_1\otimes b_1\otimes a_2\otimes b_2$ swaps the second and third tensor factors. By observing that $\sigma_{23}\circ\sigma_{12}\circ\sigma_{34}=\sigma_{13}\circ\sigma_{24}\circ\sigma_{23}$, we note that when two discrete quantum groups are cocommutative, their product defined in this way is also cocommutative. This yields

\begin{corollary}\label{productcorollary}
    $\pi_1^q(\mathbb{RP}^2)\times\pi_1^q(\mathbb{RP}^2)\not\cong\pi_1^q(\mathbb{RP}^2\times\mathbb{RP}^2)$.
\end{corollary}
\begin{rmk}
    Recall that for spaces $X$ and $Y$, $\pi_1(X\times Y)\cong\pi_1(X)\times\pi_1(Y)$. This shows that this quantum fundamental group does not obey the same property for the ordinary discrete quantum group product.
\end{rmk}

We thus have an important example where $\pi_1^q(X)$ is fully quantum. We arrive at the following question:

\begin{question}
    Can we classify which spaces admit a genuinely-quantum quantum fundamental group?
\end{question}

We answer this question in the rest of the section 
by considering the orders of $\pi_1(X)$.

\begin{rmk}\label{rmk:abelian}
 The image of a cocommutative Hopf $C^*$-algebra under a Hopf algebra morphism is cocommutative. If the finite group $\pi_1(X)$ is non-abelian, then one sees easily 
    from \Cref{generaltheorem} 
    that $\pi_1^q(X)$ genuinely quantum since $\pi_1(X)$ is a subgroup of $\pi_1^q(X)$ and 
    $C_0(\pi_1(X)) = C(\pi_1(X)) $ is non-cocommutative.     
\end{rmk}

Extending the technique of the proof of \Cref{nonclassical}, we have 

\begin{thm}\label{nonprimecase}
    If the order of $\pi_1(X)$ is finite and not prime, then $C_0(\pi_1^q(X))$ is non-cocommutative. 
\end{thm}
\begin{proof} Suppose $|\pi_1(X)|$ is not prime. If $\pi_1(X)$ is non-abelian, see \Cref{rmk:abelian}. We assume $\pi_1(X)$ is abelian. 

Pick a non-identity element $h\in\pi_1(X)$ that generates a proper (cyclic) subgroup $\Gamma<\pi_1(X)$. Note that $\Gamma$ acts on $\widetilde{X}$ by restricting the action $\alpha$ of $\pi_1(X)$ on $\widetilde{X}$. Consider the space $X':=\widetilde{X}/\Gamma$. This is the intermediate covering $p': X'\rightarrow X$ associated to the subgroup $\Gamma$. The quotient map $q:\widetilde{X}\to X'$ is a covering map, since $\widetilde{X}$ is Hausdorff and $\Gamma$ acts freely and is finite. We have $p=p'\circ q$. 

    We denote the action of an element $g\in\pi_1(X)$ via $\alpha$ by $\alpha_g:\widetilde{X}\to\widetilde{X}$ as previously. Pick an element $g\not\in\Gamma$, and denote by $\tau=(\mathbb{C},u^\tau)$ the classical deck transformation whose matrix is 
    \[u^\tau_{\tilde{x},\tilde{y}}=\delta_{\tilde{y},\alpha_g(\tilde{x})}1.\]
    Pick four points $\tilde{y}_0,\tilde{y}_1,\tilde{y}_2,\tilde{y}_3\in\tilde{X}$ which are all in distinct fibers of the covering $\tilde{X}|X$. Observe that, since $g\neq h^n$ for any $n\in\ZZ$ and the action of $\pi_1(X)$ is free, these satisfy 
    \[\{q(\tilde{y}_0),q(\tilde{y}_1),q(\tilde{y}_2),q(\tilde{y}_3)\}\cap\{q(\alpha_g(\tilde{y}_0)),q(\alpha_g(\tilde{y}_1)),q(\alpha_g(\tilde{y}_2)),q(\alpha_g(\tilde{y}_3))\}=\emptyset\]
    and $q(\tilde{y}_{i})\neq q(\tilde{y}_j)$ for any $i,j\in \{0,1,2,3\}$.

    Define a continuous, projection-valued function $P: X'\rightarrow\Bl(\CC^2)$ satisfying the following: 
    \begin{enumerate}[\rm (i)]
        \item $P(q(\alpha_g(\tilde{y}_0)))=P(q(\tilde{y}_0))=e_{11}=P(q(\tilde{y}_3))$,
        \item $P(q(\alpha_g(\tilde{y}_1)))=P(q(\tilde{y}_1))=e_{22}=P(q(\alpha_g(\tilde{y}_3)))$, and
        \item $P(q(\alpha_g(\tilde{y}_2)))=P(q(\tilde{y}_2))=t_{\frac{1}{\sqrt{2}}(e_1+e_2),\frac{1}{\sqrt{2}}(e_1+e_2)}$.
    \end{enumerate}
    This is achieved by a technique similar to that in the proof of \Cref{nonclassical}.

    In view of \Cref{formulalemmagen}, set $P_h(\tilde{x}):=P\circ q (\tilde{x}) $ and $P_e(\tilde{x})=(P_h(\tilde{x}))^\perp$ for each $\tilde{x}\in \widetilde{X}$. Note that $P_h(\tilde{x})=P_h(\alpha_{h^{-1}}(\tilde{x}))$ for all $\tilde{x}\in\widetilde{X}$, since $\Gamma$ is generated by $h$. Then the corresponding quantum deck transformation $\sigma=(\CC^2,u^{\sigma})$ is given by \[u^\sigma_{\tilde{x},\tilde{y}}:=\delta_{\tilde{y},\tilde{x}}P_e(\tilde{x})+\delta_{\tilde{y},\alpha_h(\tilde{x})}P_h(\tilde{x})\] (i.e. $P_\gamma^\sigma=0$ for all $\gamma\neq h,e$). Moreover, $\sigma$ is irreducible by \Cref{generaltheorem}. 

    Observe that \begin{align*}
        u^{\sigma\boxtimes\tau}_{\tilde{x},\tilde{y}}&=\sum_{\tilde{z}}u^\sigma_{\tilde{x},\tilde{z}}\otimes u^\tau_{\tilde{z},\tilde{y}}\\&=\sum_{\tilde{z}}\delta_{\tilde{y},\alpha_g(\tilde{z})}(\delta_{\tilde{z},\tilde{x}}P_e(\tilde{x})\otimes1+\delta_{\tilde{z},\alpha_h(\tilde{x})}P_h(\tilde{x})\otimes1)\\&=\delta_{\tilde{y},\alpha_g(\tilde{x})}P_e(\tilde{x})\otimes 1+\delta_{\tilde{y},\alpha_{gh}(\tilde{x})}P_h(\tilde{x})\otimes1
    \end{align*}
    and \begin{align*}
                u^{\tau\boxtimes\sigma}_{\tilde{x},\tilde{y}}&=\sum_{\tilde{z}}u^\tau_{\tilde{x},\tilde{z}}\otimes u^\sigma_{\tilde{z},\tilde{y}}\\&=\sum_{\tilde{z}}\delta_{\tilde{z},\alpha_g(\tilde{x})}(\delta_{\tilde{y},\tilde{z}}1\otimes P_e(\tilde{z})+\delta_{\tilde{y},\alpha_h(\tilde{z})}1\otimes P_h(\tilde{z}))\\&=\delta_{\tilde{y},\alpha_g(\tilde{x})}1\otimes P_e(\alpha_g(\tilde{x}))+\delta_{\tilde{y},\alpha_{hg}(\tilde{x})}1\otimes P_h(\alpha_g(\tilde{x})).
    \end{align*} By the usual identifications, $\sigma\boxtimes\tau\sim{\sigma'}$ and $\tau\boxtimes\sigma\sim\sigma''$ via a unitary equivalence, where ${\sigma'}=(\CC^2,u^{\sigma'})$ and $\sigma''=(\CC^2,u^{\sigma''})$ are defined via \[u^{\sigma'}_{\tilde{x},\tilde{y}}=\delta_{\tilde{y},\alpha_g(\tilde{x})}P_e(\tilde{x})+\delta_{\tilde{y},\alpha_{gh}(\tilde{x})}P_h(\tilde{x})\] and \[u^{\sigma''}_{\tilde{x},\tilde{y}}=\delta_{\tilde{y},\alpha_g(\tilde{x})}P_e(\alpha_g(\tilde{x}))+\delta_{\tilde{y},\alpha_{hg}(\tilde{x})}P_h(\alpha_g(\tilde{x})).\] 
    
    Suppose that we have a unitary equivalence $\sigma\boxtimes\tau\sim\tau\boxtimes\sigma$. Then, ${\sigma'}\sim\sigma''$, and there is a unitary $U:\CC^2\rightarrow\CC^2$ satisfying $Uu^{\sigma'}_{\tilde{x},\tilde{y}}U^*=u^{\sigma''}_{\tilde{x},\tilde{y}}$ for all $\tilde{x},\tilde{y}\in\widetilde{X}$. Fix $\tilde{x}$ and define $\tilde{y}=\alpha_{gh}(\tilde{x})$.

    Since $gh=hg$, $u^{\sigma''}_{\tilde{x},\alpha_{gh}(\tilde{x})}=P_h(\alpha_g(\tilde{x}))$. Then, \[UP_h(\tilde{x})U^*=P_h(\alpha_g(\tilde{x})).\] Since the choice of $\tilde{x}$ was arbitrary, this holds for all $\tilde{x}\in\widetilde{X}$. Choosing $\tilde{x}=\tilde{y}_0$ and $\tilde{x}=\tilde{y}_1$ in the formula above, from the construction of $P_h$, we see that $\CC\cdot e_1$ and $\CC\cdot e_2$ are invariant subspaces for $U$, so $U$ is diagonal. But then, we must also have 
    \[e_{11}=Ue_{11}U^*=UP_h(\tilde{y}_3)U^*=P_h(\alpha_g(\tilde{y}_3))=e_{22},\]
    which is impossible. 

    Thus, $\sigma\boxtimes\tau$ is not unitarily equivalent to $\tau\boxtimes\sigma$, and $C_0(\pi_1^q(X))$ is non-cocommutative.

\end{proof}

\begin{rmk}

    \begin{enumerate}[\rm (i)]
     \item The preceding proof provides an alternative to the proof of \Cref{nonclassical}. The construction is simpler and more direct here, but we preserve both proofs to provide additional concrete techniques for constructing quantum deck transformations. 
   
    \item 
    We may interpret the preceding proof as demonstrating that the rigid $C^*$-tensor categories we consider are not symmetric tensor categories.

    \end{enumerate}

\end{rmk}

What remains to be explored is the following:

\begin{question}
    Are there spaces $X$ with $\pi_1(X)\cong\ZZ/p\ZZ$, $p\geq 3$ prime and 
    $C_0(\pi_1^q(X))$
    non-cocommutative?
\end{question}

The rest of this section gives an affirmative answer.

Henceforth in this section, we suppose $\pi_1(X)\cong\ZZ/N\ZZ$ for some natural number $N\geq2$ (which includes the case where $N=p$ is prime). In this case, we have the following identification, formalizing \Cref{representations}: 

\begin{prop}\label{equivchar}
    Let $X$ be a space with fundamental group $\ZZ/N\ZZ$ and $\widetilde{X}$ be its universal cover. Let $\phi\in U_N$ be defined by $\phi_{ij}=\delta_{i,j+1}$ for $i,j\in\ZZ/N\ZZ$ be a generator of $\ZZ/N\ZZ\leq U_{N}$, and let $\alpha_j$ denote the action of $j\in\ZZ/N\ZZ$ on $\widetilde{X}$ by a deck transformation. Let $e_{ij}$, $i,j\in\{0,1,\dots, N-1\}$, denote the matrix units of $M_N(\CC)$.
    
    A quantum deck transformation $\sigma=(\Hl_\sigma,u^\sigma)$ gives rise to a family of non-degenerate representations $\rho_{\tilde{x}}^\sigma:C(S_N^+)\rightarrow\Bl(\Hl_\sigma)$ indexed by $\tilde{x}\in\widetilde{X}$ which satisfies the following: 
    \begin{enumerate}[\rm (i)]
        \item The family $\{\rho^\sigma_{\tilde{x}}\}_{\tilde{x}\in\widetilde{X}}$ is continuous in $\tilde{x}$ with respect to the topology on the set \[\{\rho:C(S_N^+)\rightarrow\Bl(\Hl_\sigma)|\text{$\,\rho$ is a $*$-homomorphism}\}\] defined via sequences by $\rho_n\rightarrow\rho$ if and only if $\rho_n(a)\rightarrow\rho(a)$ for all $a\in C(S_{N}^+)$. 
        \item For each $\tilde{x}\in\widetilde{X}$, we have 
        \begin{equation}\label{commuting}
        (\phi\otimes I_{\Hl_\sigma})(e_{ij}\otimes\rho^\sigma_{\tilde{x}}(u_{i+1,j+1}))(\phi^*\otimes I_{\Hl_\sigma})=e_{ij}\otimes\rho^\sigma_{\alpha_1(\tilde{x})}(u_{i+1,j+1})
        \end{equation} 
        for all $i,j\in\{0,1,\dots,N-1\}$.
    \end{enumerate}
    Conversely, any such family of representations satisfying these conditions gives rise to a quantum deck transformation. 
\end{prop}
\begin{proof}
    In one direction, given a quantum deck transformation $\sigma$ and $\tilde{x}\in\widetilde{X}$, we can consider the matrix \[T_{\tilde{x}}=\begin{bmatrix}
        P_{j-i}^\sigma(\alpha_{i}(\tilde{x}))
    \end{bmatrix}_{i,j=0}^{N-1}\] indexed by $\ZZ/N\ZZ$, with $i$ indexing rows and $j$ indexing columns, where the $P_k^\sigma$ are the functions in \Cref{formulalemmagen}. The conditions there show the entries of this matrix must satisfy the relations \Cref{MagicProj} and \Cref{MagicProjSum} as \Cref{matrixrmk} indicates. Then, we define a representation 
    \[\rho^\sigma_{\tilde{x}}:C(S_N^+)\rightarrow\Bl(\Hl_\sigma)\] by \[\rho^\sigma_{\tilde{x}}(u_{kl})=P_{l-k}^\sigma(\alpha_{k-1}(\tilde{x}))\] for each $u_{kl}\in C(S_N^+)$ as in \Cref{cqperm}.
    (On the right-hand side, we identify $k$ and $l$ with their equivalence classes in $\ZZ/N\ZZ$ with addition modulo $N$.)
    
    Since all $\rho_{\tilde{x}}^\sigma$ are $*$-homomorphisms, the entrywise continuity of $T_{\tilde{x}}$ in $\tilde{x}\in\widetilde{X}$ guarantees the continuity specified in the above statement. It is a straightforward computation to verify that $T_{\alpha_1(\tilde{x})}=\phi T_{\tilde{x}}\phi^*$, and the entries of $T_{\tilde{x}}$ correspond to $e_{ij}\otimes \rho^\sigma_{\tilde{x}}(u_{i+1,j+1})$ 
    via the identification $M_N(\CC)\otimes\Bl(\Hl_\sigma)\cong M_N(\Bl(\Hl_\sigma))$ for $i,j=0,\dots,N-1$.

    On the other hand, suppose we have a family of representations $\rho_{\tilde{x}}:C(S_N^+)\rightarrow\Bl(\Hl)$ satisfying the conditions above. Then, define $T_{ij}(\tilde{x}):=\rho_{\tilde{x}}(u_{i+1,j+1})$ for $i,j\in\ZZ/N\ZZ$ on the left-hand side and $i,j\in\ZZ$ on the right-hand side. Assemble the matrix \[T_{\tilde{x}}=\begin{bmatrix}
        T_{ij}(\tilde{x})
    \end{bmatrix}_{i,j=0}^{N-1},\] where $i$ labels rows and $j$ labels columns. Note that $T_{ij}(\tilde{x})$ defines a continuous, projection-valued function for each $i,j\in\ZZ/N\ZZ$. 
    
    By assumption on $\rho_{\tilde{x}}$ and the identification $M_N(\CC)\otimes\Bl(\Hl_\sigma)\cong M_N(\Bl(\Hl_\sigma))$ as $C^*$-algebras, $\phi T_{\tilde{x}}\phi^*=T_{\alpha_1(\tilde{x})}$ for all $\tilde{x}\in\widetilde{X}$. This means that \[T_{i+1,j+1}(\tilde{x})=(T_{\tilde{x}})_{i+1,j+1}=(T_{\alpha_1(\tilde{x})})_{ij}=T_{ij}(\alpha_1(\tilde{x})),\] where addition is in the group $\ZZ/N\ZZ$. 
    
    Define $P_i(\tilde{x}):=T_{0,i}(\tilde{x})$, the entry of $T_{\tilde{x}}$ in the $0^{\text{th}}$ row and $i^{\text{th}}$ column. Since the rows and columns of $T_{\tilde{x}}$ add up to the identity, by summing over the $0^{\text{th}}$ row and the $0^{\text{th}}$ column, 
    \[
    P_0(\tilde{x})+P_1(\tilde{x})+\dots+P_{N-1}(\tilde{x})=I=\sum_{j=0}^{N-1}P_{-j}(\alpha_{j}(\tilde{x}))
    \]
    in view of the above. This implies that the functions $P_i=P_i^\sigma$ for some quantum deck transformation $\sigma$, whose underlying Hilbert space is $\Hl_\sigma:=\Hl$, by \Cref{formulalemmagen}. 
\end{proof}

We will continue to use a superscript to denote the family 
$\{\rho^\sigma_{\tilde{x}}\}_{\tilde{x}\in\widetilde{X}}$ associated to a quantum deck transformation $\sigma$ as above.

\begin{lemma}\label{periodicity}
    $C_0(\pi_1^q(X))$ is cocommutative if and only if, for every irreducible quantum deck transformation $\sigma$, the family $\{\rho_{\tilde{x}}^\sigma\}$ of \Cref{equivchar} 
    satisfies $\rho_{\alpha_j(\tilde{x})}^\sigma=\rho_{\tilde{x}}^\sigma$ for 
    all $\tilde{x} \in \widetilde{X} $ and $j \in \ZZ/N\ZZ$. 
\end{lemma}
\begin{proof}
    If for each irreducible quantum deck transformation $\sigma_0$, $\rho_{\alpha_j(\tilde{x})}^{\sigma_0}=\rho_{\tilde{x}}^{\sigma_0}$ for all $j\in\ZZ/N\ZZ$, it is a routine computation (as in the $\mathbb{RP}^2$ case) to see that $\sigma\boxtimes\tau\sim\tau\boxtimes\sigma$ by conjugation by $t$ for all irreducible quantum deck transformations $\sigma$ and $\tau$, implying $\phi\circ\Delta=\Delta$.
    
    Conversely, suppose $C_0(\pi_1^q(X))$ is cocommutative. This implies that $t$ implements an equivalence $\sigma\boxtimes\tau\sim\tau\boxtimes\sigma$ for all irreducible quantum deck transformations $\sigma$ and $\tau$. 

    Let $\sigma$ be arbitrary and irreducible, and set $\tau=(\CC,u^\tau)$, where $u^\tau_{\tilde{x},\tilde{y}}:=\delta_{\tilde{y},\alpha_1(\tilde{x})}$ for each $\tilde{x},\tilde{y}\in\widetilde{X}$. Then, writing $u^\sigma_{\tilde{x}',\tilde{y}'}=\sum_{k=0}^{N-1}\delta_{\tilde{y}',\alpha_k(\tilde{x}')}P_k^\sigma(\tilde{x}')$ for each $\tilde{x}',\tilde{y}'\in\widetilde{X}$, we have 
    
    \begin{align*}
u^{\sigma\boxtimes\tau}_{\tilde{x},\tilde{y}}
&=\sum_{\tilde{z}\in\widetilde{X}}u^\sigma_{\tilde{x},\tilde{z}}\otimes u^\tau_{\tilde{z},\tilde{y}}
    =\sum_{\tilde{z}\in\widetilde{X}}\sum_{k=0}^{N-1}\delta_{\tilde{z},\alpha_k(\tilde{x})}\delta_{\tilde{y},\alpha_1(\tilde{z})}P_k^\sigma(\tilde{x})\otimes1\\
    &=\sum_{k=0}^{N-1}\delta_{\tilde{y},\alpha_{k+1}(\tilde{x})}P_k^\sigma(\tilde{x})\otimes 1
    \end{align*} 
    and 
    \begin{align*}
    u^{\tau\boxtimes\sigma}_{\tilde{x},\tilde{y}}
    &=\sum_{\tilde{z}\in\widetilde{X}}u^\tau_{\tilde{x},\tilde{z}}\otimes u^\sigma_{\tilde{z},\tilde{y}}
    =\sum_{\tilde{z}\in\widetilde{X}}\sum_{k=0}^{N-1}\delta_{\tilde{z},\alpha_1(\tilde{x})}\delta_{\tilde{y},\alpha_k(\tilde{z})}1\otimes P_k^\sigma(\tilde{z})\\
    &=\sum_{k=0}^{N-1}\delta_{\tilde{y},\alpha_{k+1}(\tilde{x})}1\otimes P_k^\sigma(\alpha_1(\tilde{x}))
    \end{align*} 
    for all $\tilde{x},\tilde{y}\in\widetilde{X}$. 
    Since $tu_{\tilde{x}',\tilde{y}'}^{\tau\boxtimes\sigma}t=u^{\sigma\boxtimes\tau}_{\tilde{x}',\tilde{y}'}$ for all $\tilde{x}',\tilde{y}'\in\widetilde{X}$, fixing $\tilde{x}$ and $k\in\ZZ/N\ZZ$ and setting $\tilde{y}=\alpha_{k+1}(\tilde{x})$,  we deduce \[P_k(\alpha_1(\tilde{x}))=P_k(\tilde{x}).\] Since $\tilde{x}$ and $k$ are arbitrary and $\alpha_1$ generates the deck transformation group of $\tilde{X}$, we obtain that $P_k\circ\alpha_j=P_k$ for all $k,j\in\ZZ/N\ZZ$. 
    This implies that $\rho_{\alpha_j(\tilde{x})}^\sigma=\rho_{\tilde{x}}^\sigma$ for all $\tilde{x}\in\widetilde{X}$. 
\end{proof}
\begin{rmk}
This provides another route to prove \Cref{nonprimecase}. 
\end{rmk}

In the first version of this paper, 
Theorem 6.14 mistakenly stated that 
$C_0(\pi_1^q(X))$ is cocommutative if $\pi_1(X)\cong\ZZ/3\ZZ$. 
One of the authors miraculously received without solicitation in his ChatGPT 
account a counter-example to the original Theorem 6.14 that was posted on the arXiv. 
We generalized the construction in the counter-example to reach 
the current version of \Cref{thm:qfgpnccm}.

\begin{thm} \label{thm:qfgpnccm}
    For each odd prime $p$, there exists a space $X_p$ such that $\pi_1(X_p)\cong\ZZ/p\ZZ$ and $C_0(\pi_1^q(X_p))$ is non-cocommutative.
\end{thm}
\begin{proof}
    Let $p$ be an odd prime, and let $\omega:=e^{2\pi i/p}$. Define \[\gamma:=\begin{bmatrix}
        \omega & 0 & \cdots & 0 \\
        0 & \omega^2 & \cdots & 0 \\
        \vdots & \vdots &\ddots & \vdots \\
        0 & 0 & \cdots & \omega^{p-1}
    \end{bmatrix}\in SU_{p-1},\] and let $\Gamma$ denote the cyclic subgroup of $SU_{p-1}$ generated by $\gamma$. Then, $\Gamma$ acts on $SU_{p-1}$ by right action via matrix multiplication. Denote this action by $\alpha_{g}$ for each $g\in\Gamma\cong\ZZ/p\ZZ=\{0,1,\dots,p-1\}$. 
    
    Let $X_p:=SU_{p-1}/\Gamma$. Since $SU_{p-1}$ is compact Hausdorff and $\Gamma$ acts freely, the quotient map $q:SU_{p-1}\to X_p$ is a covering map. Moreover, $SU_{p-1}$ is simply connected, so it is the universal cover for $X_p$. Hence, $\pi_1(X_p)=\Gamma\cong\ZZ/p\ZZ$.  

    Let \[f_{j}:=\frac{1}{\sqrt{p}}(1,\omega^j,\omega^{2j},\dots,\omega^{(p-1)j})\] 
    viewed as column matrices for $j\in\ZZ/p\ZZ$. The set $\{f_j\}_{j\in\ZZ/p\ZZ}$ is an orthonormal basis for $\CC^p$. Let $Q_j:=t_{f_j,f_j}$ denote the associated collection of pairwise orthogonal projections, which sum to the identity. Define the map $\psi:SU_{p-1}\to SU_{p}$ by $\psi(\tilde{x})=1\oplus\tilde{x}$. Since $\psi(\gamma)f_j=(1\oplus\gamma)f_j=f_{j+1}$, 
    \begin{equation}\label{eqn:cyclicvec}
    \psi(\gamma)Q_j\psi(\gamma)^*=t_{\psi(\gamma)f_j,\psi(\gamma)f_j}=t_{f_{j+1},f_{j+1}}=Q_{j+1}
    \end{equation}
    for all $j\in\ZZ/p\ZZ$, addition modulo $p$.

    For each $\tilde{x}\in SU_{p-1}$ and $k\in \ZZ/p\ZZ$, let 
    \begin{equation}\label{eqn:Pkdefn}
    P_k(\tilde{x}):=\psi(\tilde{x})Q_{2k}\psi(\tilde{x})^*,
    \end{equation} where the multiplication $2k$ is modulo $p$. This defines a collection $\{P_k\}_{k\in\ZZ/p\ZZ}$ of continuous, projection-valued functions on $SU_{p-1}$. 
    
    Note that $k\mapsto2k\mod p$ is a permutation of $\ZZ/p\ZZ$, since $p$ is prime. Then, for each $\tilde{x}$, 
    \begin{equation}\label{condition1}
    P_0(\tilde{x})+P_1(\tilde{x})+\dots+P_{p-1}(\tilde{x})=\psi(\tilde{x})\bigg(\sum_{k=0}^{p-1}Q_{2k}\bigg)\psi(\tilde{x})^*=\psi(\tilde{x})\bigg(\sum_{k=0}^{p-1}Q_{k}\bigg)\psi(\tilde{x})^*=I_{p\times p}.
    \end{equation} For each $k\neq0$, by repeatedly applying \Cref{eqn:cyclicvec} and using that $\psi$ is a group homomorphism,
    \begin{align}
    P_k(\alpha_{p-k}(\tilde{x}))
    &=P_k(\tilde{x}\gamma^{p-k})=\psi(\tilde{x})\psi(\gamma)^{p-k}Q_{2k}\psi(\gamma)^{*p-k}\psi(\tilde{x})^*\notag\\
    &=\psi(\tilde{x})Q_{2k+p-k}\psi(\tilde{x})^*=\psi(\tilde{x})Q_k\psi(\tilde{x})^*,\label{eqn:composition}
    \end{align}
    since $p-k+2k\equiv k\mod p$. Since $k\mapsto2k\mod p$ maps $\ZZ/p\ZZ-\{0\}$ to itself, we have 
    \begin{equation}\label{eqn:seteq1}
    \bigg\{\psi(\tilde{x})Q_{k}\psi(\tilde{x})^*:k=1,\dots,p-1\bigg\}= \bigg\{\psi(\tilde{x})Q_{2k}\psi(\tilde{x})^*:k=1,\dots,p-1\bigg\}
    \end{equation} for all $\tilde{x}\in SU_{p-1}$. By \Cref{eqn:Pkdefn}, \Cref{eqn:composition}, and \Cref{eqn:seteq1}, the sets 
    \[
   \bigg \{P_k(\alpha_{p-k}(\tilde{x})):k=1,\dots,p-1\bigg\}=\bigg\{P_k(\tilde{x}):k=1,\dots,p-1\bigg\}
    \] are equal for all $\tilde{x}\in SU_{p-1}$. It follows that 
    \begin{equation}\label{condition2}
    \sum_{k=1}^{p-1}P_k(\alpha_{p-k}(\tilde{x}))=\sum_{k=1}^{p-1}P_k(\tilde{x})
    \end{equation} for all $\tilde{x}\in SU_{p-1}$. By \Cref{formulalemmagen}, \Cref{condition1}, and \Cref{condition2}, we obtain a quantum deck transformation $\sigma$. 

    We now check that this $\sigma$ is irreducible. By \Cref{generaltheorem}, it suffices to show that \[A:=C^*(P_k(\tilde{x}):\tilde{x}\in SU_{p-1},k\in\ZZ/p\ZZ)=M_{p}(\CC).\] Let $2^{-1}:=(p+1)/2$ denote the multiplicative inverse of $2$ in the field $\ZZ/p\ZZ$. We can rearrange \Cref{eqn:Pkdefn} to \begin{equation}
        \label{eqn:Pkdefnalt}
        \psi(\tilde{x})Q_j\psi(\tilde{x})^*=P_{2^{-1}j}(\tilde{x})
    \end{equation}
    for all $j\in\ZZ/p\ZZ$.
    
    Fix 
    \[
    \tilde{x}_0:=\begin{bmatrix}
        -\omega^{p-1} & 0 & \cdots & 0 & 0\\
        0 & -\omega^{p-2} & \cdots & 0 & 0 \\
        \vdots & \vdots & \ddots & \vdots & \vdots \\
        0 & 0 & \cdots & -\omega^2 & 0 \\
        0 & 0 & \cdots & 0 & -\omega
    \end{bmatrix}\in SU_{p-1}.
    \]
    Note that
    \[
  v := \psi(\tilde{x}_0)f_1=\frac{1}{\sqrt{p}}(1,-1,-1,\dots,-1,-1)\in\CC^p.
    \] Let $Q:=t_{v,v}$ denote the projection to the line $\CC v\subset\CC^p$. Since $Q=\psi(\tilde{x}_0)Q_1\psi(\tilde{x}_0)^*=P_{2^{-1}}(\tilde{x}_0)$ by \Cref{eqn:Pkdefnalt}, we have $Q\in A$. 

    By \Cref{eqn:Pkdefnalt} again, $Q_j=P_{2^{-1}j}(I_{p-1\times p-1})\in A$ for all $j\in\ZZ/p\ZZ$. Let $k,l,m\in\ZZ/p\ZZ$. We have
    \[
    Q_kQQ_lf_m=\delta_{lm}Q_kQf_l=\delta_{lm}(f_l,v)Q_kv=\delta_{lm}(f_l,v)(v,f_k)f_k=\delta_{lm}(f_l,v)\overline{(f_k,v)}f_k.
    \] 
    For any $j\in\ZZ/p\ZZ$,
    \[
        (f_j,v)=\frac{1}{p}\left(1+\sum_{k=1}^{p-1}-\omega^{kj}\right)=\frac{2}{p}\neq0.
    \]
    Thus, $(f_l,v)(v,f_k)\neq0$. Since $Q_kQQ_l\in A$, we obtain that the matrix units $t_{f_k,f_l}\in A$ for all $k,l\in\ZZ/p\ZZ$. Hence, $A=M_p(\CC)$, and $\sigma$ is irreducible. 

    Now, by \Cref{periodicity}, it suffices to show that the family of representations $\{\rho^\sigma_{\tilde{x}}\}_{\tilde{x}\in SU_{p-1}}$ is not invariant under the action $\alpha$. We have that $\rho_{\tilde{x}}^\sigma(u_{11})=P_0(\tilde{x})$, while \[
    \rho_{\alpha_1(\tilde{x})}^\sigma(u_{11})=P_0(\alpha_1(\tilde{x}))=\psi(\tilde{x})\psi(\gamma)Q_0\psi(\gamma)^*\psi(\tilde{x})^*=\psi(\tilde{x})Q_1\psi(\tilde{x})^*\neq\rho^\sigma_{\tilde{x}}(u_{11}).\]
\end{proof}

\section{Galois Correspondence}\label{sec:galois}

Here, we address the second question at the beginning of \Cref{sec:clas}: which properties of the classical fundamental group carry over to the quantum setting? \Cref{productcorollary} indicates there are important differences. Nevertheless, we conclude with a Galois-type correspondence.

Throughout, $X$ is a locally compact Hausdorff space which is path-connected, locally path-connected, and semi-locally simply connected. Also, let $p:\widetilde{X}\rightarrow X$ be a universal covering map of $X$. We do not require $|\pi_1(X)|<\infty$.

The classical Galois correspondence holds in one direction:

\begin{thm}\label{galois}
    Let $p_1:X_1\rightarrow X$ be a covering space, and let $q:\widetilde{X}_1\rightarrow X_1$ be the universal cover of $X_1$. Then, there exists a surjective $*$-homomorphism $\phi:C_0(\pi_1^q(X))\twoheadrightarrow C_0(\pi_1^q(X_1))$ which is compatible with the comultiplications. 
    That is, $\pi_1^q(X_1)$ is a discrete quantum subgroup of $\pi_1^q(X)$.
\end{thm}

\begin{proof}
    The argument here is a modification of the one to show $\pi_1^q(X)$ is well-defined in \Cref{welldefined}, with $F$ there replaced by $ r$ to be defined below.  
    
    Since $X$ satisfies the aforementioned properties, the composition $p_1\circ q:\widetilde{X}_1\rightarrow X$ is a covering map. Since $\widetilde{X}_1$ is simply connected, we have a homeomorphism (in fact, an isomorphism of covering spaces) $r:\widetilde{X}_1\rightarrow\widetilde{X}$. Thus, we have the following commutative diagram: 
    \[\begin{tikzcd}
        \widetilde{X}_1 \arrow{r}{q}  \arrow{d}{r} & X_1 \arrow{d}{p_1^{}} \\
        \widetilde{X} \arrow{r}{p} & X
    \end{tikzcd}
    \]

    Let $\sigma=(\Hl_\sigma,u^\sigma)\in\text{Deck}^+(\tilde{X_1}|X_1)$. We define $r_*(\sigma):=(\Hl_\sigma,v^{r_*(\sigma)})$ by setting 
    \[v^{ r_*(\sigma)}_{\tilde{x},\tilde{y}}:=u^\sigma_{ r^{-1}(\tilde{x}), r^{-1}(\tilde{y})}\] 
    for each $\tilde{x},\tilde{y}\in\widetilde{X}$. We have that, since $ r$ is surjective, the matrix $v^{ r_*(\sigma)}$ satisfies the relations to be a quantum permutation. Moreover, since $ r$ and $ r^{-1}$ are continuous, $ r_*(\sigma)$ is a quantum homeomorphism.

    Suppose $p(\tilde{x})\neq p(\tilde{y})$. Then,
    $p_1(q\circ  r^{-1}(\tilde{x}))\neq p_1(q\circ  r^{-1}(\tilde{y}))$, 
    hence $q( r^{-1}(\tilde{x}))\neq q( r^{-1}(\tilde{y}))$. It follows that $u^{\sigma}_{ r^{-1}(\tilde{x}), r^{-1}(\tilde{y})}=0$, 
    since $\sigma\in\text{Deck}^+(\widetilde{X}_1|X_1)$. Therefore, $ r_*(\sigma)$ is an object in $\text{Deck}^+(\widetilde{X}|X)$. Since $ r$ is bijective, if we define $ r_*$ to be the identity on morphisms, \begin{equation}\label{identitymorphisms}
    \Mor( r_*(\sigma), r_*(\tau))=\Mor(\sigma,\tau),
    \end{equation} 
     we obtain a functor
    \[
     r_*:\text{Deck}^+(\widetilde{X}_1|X_1)\rightarrow\text{Deck}^+(\widetilde{X}|X).
    \]

    By \eqref{identitymorphisms}, we have that $ r_*$ maps irreducible objects to irreducible objects and inequivalent objects to inequivalent objects. That is, $ r_*$ defines an injective map on unitary equivalence classes of irreducible objects. Moreover, $ r_*$ is also compatible with the tensor structure, since for all $\tilde{x},\tilde{y}\in\widetilde{X}$, 
    \begin{align*}
    v^{ r_*(\sigma\boxtimes\tau)}_{\tilde{x},\tilde{y}}
    &=u^{\sigma\boxtimes\tau}_{ r^{-1}(\tilde{x}), r^{-1}(\tilde{y})}
    =\sum_{\tilde{z}_1\in\widetilde{X}_1}u^\sigma_{ r^{-1}(\tilde{x}),\tilde{z}_1}\otimes u^\tau_{\tilde{z}_1, r^{-1}(\tilde{y})}\\
    &=\sum_{\tilde{z}\in\widetilde{X}}u^\sigma_{ r^{-1}(\tilde{x}), r^{-1}(\tilde{z})}\otimes u^\tau_{ r^{-1}(\tilde{z}), r^{-1}(\tilde{y})}
    =v^{ r_*(\sigma)\boxtimes \:r_*(\tau)}_{\tilde{x},\tilde{y}}.
    \end{align*}
    Let $ r_*(\text{Irred}_1^q(X_1))\subseteq\text{Irred}_1^q(X)$ denote the image of the elements of $\text{Irred}_1^q(X_1)$ under $ r_*$. Then, 
    \begin{align*}
v_{\tilde{x},\tilde{y}}&=\sum_{\sigma'\in  r_*(\text{Irred}_1^q(X_1))}v^{\sigma'}_{\tilde{x},\tilde{y}}+\sum_{\sigma'\not\in  r_*(\text{Irred}_1^q(X_1))}v^{\sigma'}_{\tilde{x},\tilde{y}}\\
&=\sum_{\sigma\in\text{Irred}_1^q(X_1)}v^{  r_*(\sigma)}_{\tilde{x},\tilde{y}}+\sum_{\sigma'\not\in  r_*(\text{Irred}_1^q(X_1))}v^{\sigma'}_{\tilde{x},\tilde{y}}.
     \end{align*}
    Let $e_{ij}^{\sigma'}$ denote the matrix units of $C_0(\pi_1^q(X))$ and $f_{ij}^\sigma$ denote the matrix units of $C_0(\pi_1^q(X_1))$. Define $\phi^*:C_0(\pi_1^q(X))\rightarrow C_0(\pi_1^q(X_1))$ by 
    \[    \phi^*(e_{ij}^{\sigma'})=0
    \] 
    for $\sigma'\not\in  r_*(\text{Irred}_1^q(X_1))$ and
    \[\phi^*(e^{\sigma'}_{ij})=f^\sigma_{ij}
    \] 
    for $\sigma'=  r_*(\sigma)\in  r_*(\text{Irred}_1^q(X_1))$. By the injectivity of $  r_*:\text{Irred}_1^q(X_1)\to \text{Irred}_1^q(X)$, this determines a surjective $*$-homomorphism $\phi^*:C_0(\pi_1^q(X))\rightarrow C_0(\pi_1^q(X_1))$ satisfying $\phi^*(v_{\tilde{x},\tilde{y}})=u_{\tilde{x},\tilde{y}}$. Moreover, since $  r_*$ preserves the tensor structure and irreducibility, and \[
    u_{\tilde{x},\tilde{y}}=\sum_{\sigma\in\text{Irred}_1^q(X_1)}u_{\tilde{x},\tilde{y}}^\sigma,
    \] 
    the map $\phi^*$ is compatible with the comultiplications. Hence, $\pi_1^q(X_1)\leq\pi_1^q(X)$.
\end{proof}

\begin{rmk}
\label{incompleteness}
    Note that this correspondence is not as complete as in the classical situation. There, covering spaces are in one-to-one correspondence with subgroups of the fundamental group; however, it may be possible many discrete quantum subgroups of the quantum fundamental group may correspond to the same covering by an action.

    Motivated by this incomplete correspondence, we may ask: \begin{question} What is the appropriate quantum analogue of covering spaces which completes the correspondence? \end{question} That is, given a $C^*$-algebra $A$, which we regard as an algebra of functions $C_0(X)$ on a quantum space, is there some appropriate category whose objects are pairs $(B,p^*)$, where $B$ is a $C^*$-algebra and $p^*\in\Mor(A,B)$ is a $*$-homomorphism, which includes the classical notion of a covering space? Our work on this question is ongoing.
\end{rmk}

\noindent
{\bf Acknowledgments}: 
M.Z. is grateful to the Franklin College 
at The University of Georgia for a research assistantship in summer of 2026. The authors would like to thank the algebra group of our mathematics department for inviting M.Z. to give a talk in their seminar 
on our work in progress. Questions raised by Nate Harman and Daniel Nakano during the seminar led us to 
consider questions on cocommutativity in \Cref{sec:clas}. 
Gemini and ChatGPT were used for some basic concepts in algebraic topology  
and literature search. 
After the paper was submitted for publication and posted on the arXiv with the original Theorem 6.14 erroneously stating that 
$C_0(\pi_1^q(X))$ is cocommutative if $\pi_1(X)\cong\ZZ/3\ZZ$,  
one of the authors miraculously received without solicitation in his ChatGPT 
account a counter-example to the original Theorem 6.14. 
We generalized the construction in the counter-example to reach 
the current version of \Cref{thm:qfgpnccm}. 

\bigskip 
\noindent 
{\bf Declaration of competing interest:}\
The authors declare that they have no known competing financial interests or personal relationships that could have appeared to influence the work reported in this paper.

\bigskip 
\noindent 
{\bf Data availability:} \
No data were generated for this study.


\begin{thebibliography}{99}


\bibitem{Baaj-Skan93} Baaj, S. and Skandalis, G.,
Unitaires multiplicatifs et dualit\'e pour les produits crois\'es de $C^*$-alg\`ebres,
{\it Ann. Sci. Ec. Norm. Sup.} {\bf 26} (1993), no. 4, 425--488.


\bibitem{Ban05a} Banica, T.,
Quantum automorphism groups of homogeneous graphs, 
{\it J. Funct. Anal.} {\bf 224} (2005), no. 2, 243--280.


\bibitem{Ban05b}
Banica, T.,
Quantum automorphism groups of small metric spaces,
{\it Pacific J. Math.} {\bf 219} (2005), no. 1, 27--51.

\bibitem{Ban-Sp09}
Banica, T. and Speicher, R.,
Liberation of orthogonal Lie groups,
{\it Adv. Math.} {\bf 222} (2009), no. 4, 1461--1501.

\bibitem{Bich03a} Bichon, J.,
Quantum automorphism groups of finite graphs,
{\it Proc. Amer. Math. Soc.} {\bf 131} (2003), no. 3, 665--673

\bibitem{Dr86} Drinfeld, V.~G.,
{\rm Quantum groups,}
{\it Proc. ICM-1986},  Vol. 1, 2 (Berkeley, Calif., 1986), 798--820,
Amer. Math. Soc., Providence, R.I., 1987.

\bibitem{Daws12} Daws, M., Kasprzak, P., Skalski, A., So{\l}tan, P.~M.,
{\rm Closed quantum subgroups of locally compact quantum groups},
{\it Adv. Math.} {\bf 231} (2012), no. 6, 3473--3501.

\bibitem{EffrosRuan94} Effros, E. and Ruan, Z.-J.,
{\rm Discrete quantum groups, I. The Haar measure,}
{\it Internat. J. Math.} {\bf 5} (1994), no. 5, 681--723.


\bibitem{Jimb85} Jimbo, M.,
{A $q$-difference analogue of $U\mathfrak{g}$ and the Yang-Baxter equations,}  
{\it Lett.  Math. Phys.} {\bf 10}, (1985), no. 1, 63--69.

\bibitem{KusVaes00}  Kustermans, J. and Vaes, S.,
{\rm Locally compact quantum groups},
{\it Ann. Sci. {\'E}cole Norm. Sup. } (3) {\bf 33 } (2000), no. 6, 837–934.

\bibitem{KusVaes03}  Kustermans, J. and Vaes, S.,
{\rm Locally compact quantum groups in the von Neumann algebraic setting},
{\it Math. Scand.} (3) {\bf 92 } (2003), no. 1, 68--92.


\bibitem{Meyer12}
Meyer, R., Roy, S., and Woronowicz, S. L.,
{\rm  Homomorphisms of quantum groups},
{\it M{\"u}nster J. Math.} {\bf 5} (2012), 1--24.


\bibitem{Nesh_Tuset13}
Neshveyev, S. and Tuset, L., 
\textit{Compact quantum groups and their representation categories}, 
Cours Sp\'{e}cialis\'{e}s [Specialized Courses] 20, Soci\'{e}t\'{e} Math\'{e}matique de France, Paris, 2013, vi+169 pp.



\bibitem{PoWo90}
Podle\'s, P. and Woronowicz, S.~L., 
Quantum deformation of {L}orentz group,
{\it Comm. Math. Phys.}, {\bf 130} (1990), no. 2, 381--431.

\bibitem{Vaes05}
Vaes, S.,
{\rm  A new approach to induction and imprimitivity results},
\textit{J. Funct. Anal. } {\bf 229} (2005), no. 2, 317--374.


\bibitem{Vdl-Haar} 
Van Daele, A.,
{The Haar measure of a compact quantum group},
Proc. Amer. Math. Soc.
{\bf 123} (1995), 3125-3128.

\bibitem{Vdl-discr96}
Van~Daele, A., 
Discrete quantum groups,
{\it J. Algebra}, {\bf 180} (1996), no. 2, 431--444.

\bibitem{Vdl-duality98}
Van~Daele, A., 
An algebraic framework for group duality,
{\it Adv. Math.}, {\bf 140} (1998), no. 2, 323--366.

\bibitem{Vdl-discr25}
Van~Daele, A., 
Discrete quantum groups and their duals,  
 \url{https://doi.org/10.48550/arXiv.2512.12350}, (2026).

\bibitem{VanWan96} Van Daele, A. and Wang, S. Z.,
Universal quantum groups,
{\it International J. of Math.} {\bf 7} (1996), no. 2, 255--263.

\bibitem{Voigt23} 
Voigt, C., 
Infinite quantum permutations,  
\textit{ Adv. Math.} \textbf{415} (2023), Paper No. 108887,  37 pp. 

\bibitem{Voigt26} 
Voigt, C.,
An introduction to quantum symmetries, 
\url{https://doi.org/10.48550/arXiv.2603.25141}, (2026).


\bibitem{Wan95}
Wang, S. Z.,
Free products of compact quantum groups,
{\it Comm. Math. Phys.}, {\bf 167} (1995), no. 3, 671--692.

\bibitem{Wan98}
Wang, S. Z.,
Quantum symmetry groups of finite spaces, 
\text{\it Comm. Math. Phys.} \textbf{195} (1998), no. 1, 195--211.

\bibitem{Wan02b} Wang, S.~Z.,
Structure and isomorphism classification of quantum
groups $A_u(Q)$ and $B_u(Q)$,
{\it Journal of Operator Theory} {\bf 48} (2002), no. 3, suppl., 573–583.


\bibitem{Wan13}
Wang, S. Z., 
Equivalent notions of normal quantum subgroups, compact quantum
groups with properties {$F$} and {$FD$}, and other applications, 
{\it J. Algebra}, {\bf 397} (2014), 515--534.

\bibitem{Wor87-SU2} 
Woronowicz, S.~L.,  
{Twisted $SU(2)$ group. An example of noncommutative differential
	calculus},
{\it  Publ. RIMS, Kyoto Univ.} {\bf 23} (1987), no. 1, 117--181.

\bibitem{Wor87}
Woronowicz, S.~L.,
Compact matrix pseudogroups,
{\it Comm. Math. Phys.}, {\bf 111} (1987), no. 4, 613--665.


\bibitem{Wor88} 
Woronowicz, S.~L.,
Tannaka-Kre\u{i}n duality for compact matrix pseudogroups. Twisted $SU(N)$  groups,
\textit{Invent. Math.} {\bf 93} (1988), no. 1, 35--76.

\bibitem{Wor98}
Woronowicz, S.~L.,
Compact quantum groups,
In {\it Sym\'etries quantiques ({L}es {H}ouches, 1995)}, 
845--884. North-Holland, Amsterdam, 1998.



\end{thebibliography}
\end{document}